\documentclass[11pt]{amsart}

\usepackage{amsthm,amsfonts,amssymb,amsmath}
\usepackage{caption,graphicx}
\usepackage{paralist}

\usepackage{fullpage}
\usepackage{hyperref}
\usepackage{amsthm,amsfonts,amssymb,amsmath,oldgerm}
\usepackage[T1]{fontenc}
\usepackage{lmodern} 
\usepackage[francais,english]{babel}
\usepackage{amsmath,amsthm}
\usepackage{mathrsfs}
\usepackage{bbm} 
\usepackage{amssymb}
\usepackage{amsfonts}
\usepackage{pifont}
\usepackage{stmaryrd}
\usepackage{dsfont}
\usepackage{graphicx,pstricks}
\usepackage{boites,boites_exemples}
\usepackage{enumerate}
\usepackage{lscape}
\usepackage[all]{xy}
\input xy
\xyoption{all}

\usepackage{hyperref}

\newcommand{\be}{\begin{enumerate}}
\newcommand{\ee}{\end{enumerate}}
\newcommand{\pa}{\partial}

\numberwithin{equation}{section}

\newcommand{\R}{\mathbb{R}}

\newcommand{\N}{\mathbb{N}}

\newcommand{\T}{\mathbb{T}}

\newcommand{\na}{\nabla}

\newtheorem{thm}{Theorem}[section]

\newtheorem{coro}{Corollary}[section]
\newtheorem{prop}{Proposition}[section]
\newtheorem{lem}{Lemma}[section]
\newtheorem{rem}{Remark}[section]

\newcommand{\eps}{\varepsilon}

\begin{document}

\title{On the linearized Vlasov-Poisson system on the whole space around stable homogeneous equilibria}
\author{D. Han-Kwan} 
\address{Centre de Math\'ematiques Laurent Schwartz (UMR 7640), Ecole Polytechnique, Institut Polytechnique de Paris, 91128 Palaiseau Cedex, France}
\email{daniel.han-kwan@polytechnique.edu}
\author{T.T. Nguyen}
\address{Department of Mathematics, Penn State University, State College, PA 16803, USA.}
\email{nguyen@math.psu.edu}
\author{F. Rousset}
\address{Laboratoire de Math\'ematiques d'Orsay (UMR 8628), Universit\'e Paris-Saclay, 91405 Orsay Cedex, France.}
\email{frederic.rousset@universite-paris-saclay.fr}

\maketitle 

\begin{abstract}
We study the linearized Vlasov-Poisson system around suitably stable  homogeneous equilibria on $\R^d\times \R^d$ (for any $d \geq 1$) and establish  dispersive $L^\infty$ decay estimates in the physical space.
\end{abstract}

\section{Introduction}

This work is concerned with the Vlasov-Poisson system on $\R^d \times \R^d$
 for $d \geq 1$:
\begin{equation}
\label{VP}
\left \{ 
\begin{aligned}
&\pa_t \mathrm{f} + v\cdot \na_x \mathrm{f} +  \mathrm{E} \cdot \na_v  \mathrm{f}  = 0, \quad (x,v) \in \R^d \times \R^d, \\
&\mathrm{E} = \nabla_x   \Delta_x^{-1}\left( \rho- 1 \right), \qquad \rho(t,x)= \int_{\R^d} \mathrm{f} (t,x,v)\, dv,\\
&\mathrm{f}|_{t=0} = \mathrm{f}_0, %\qquad \int_{\R^d\times \R^d} \mathrm{f}_0  \, dv \, dx =1,
\end{aligned}
\right.
\end{equation}
where $\mathrm{f}$ (resp. $\mathrm{E}$) describes the distribution function of negatively charged particles (resp. the electric field) in a plasma with a fixed uniform background of ions.
We are interested in the long time behavior of the solutions to~\eqref{VP} around   \emph{homogeneous} equilibria, i.e. non-negative distribution functions $\mu(v)$ satisfying 
\begin{equation}
\label{eq:norma}
\int_{\R^d} \mu(v) \, dv = 1.
\end{equation}
To this end, we consider solutions of the form $\mathrm{f}(t,x,v) = \mu(v) + f(t,x,v)$ %, so that $f$ satisfies
%\begin{equation}
%\label{VP2}
%\left \{ 
%\begin{aligned}
%&\pa_t f + v\cdot \na_x f +  {E} \cdot \na_v (\mu +f) = 0, \quad (x,v) \in \R^d \times \R^d, \\
%&{E} = \nabla_x   \Delta_x^{-1}\rho, \qquad \rho(t,x)= \int_{\R^d} f (t,x,v)\, dv,\\
%&f|_{t=0} = f_0, \qquad \int_{\R^d\times \R^d} f_0  \, dv \, dx =0.
%\end{aligned}
%\right.
%\end{equation}
and specifically focus on the \emph{linearized} equations: %Vlasov-Poisson system  around such $\mu$, that is to say the system
\begin{equation}
\label{linVP}
\left \{ 
\begin{aligned}
&\pa_t f + v\cdot \na_x f +  E \cdot \na_v \mu = 0, \quad (x,v) \in \R^d \times \R^d, \\
&E = \nabla_x   \Delta_x^{-1} \rho, \qquad \rho(t,x)= \int_{\R^d} f (t,x,v)\, dv,\\
&f|_{t=0} = f_0.
\end{aligned}
\right.
\end{equation}
Our goal is to establish decay in time for the density $\rho$ of the solution to \eqref{linVP}. To this purpose, we will require that $\mu$ satisfies some appropriate  conditions of stability. This problem can be seen as a first step towards the understanding of relaxation properties around stable homogeneous equilibria  (i.e. \emph{Landau Damping})  for the full Vlasov-Poisson system~\eqref{VP} on the whole space.

Landau Damping was studied in the breakthrough paper \cite{MV} by Mouhot and Villani in the case of $\T^d \times \R^d$ (see also \cite{BMM1} and very recently \cite{GNR}). All these works are based on a linear mechanism called \emph{phase mixing}, which is specific to the free transport operator $
\pa_t + v \cdot \na_x
$ on the torus; furthermore they require perturbations of Gevrey or analytic regularity to handle the non-linear problem, in order to avoid resonances referred to as \emph{plasma echoes}.
For what concerns the whole space,  an important contribution is due to Bedrossian, Masmoudi and Mouhot who  considered in \cite{BMM2} the \emph{screened} Vlasov-Poisson system, which corresponds to a low frequency (or equivalently, long range) regularization of the Coulomb potential, resulting in the equation
$$
E = \nabla_x   (1-\Delta_x)^{-1} \rho
$$
for the electric field. They relied on \emph{dispersive} properties of the free transport operator
$
\pa_t + v \cdot \na_x
$
on the whole space in the Fourier side 
 to prove decay in finite regularity for the full non-linear system in dimensions $d\geq 3$ (with a strategy inspired by \cite{MV,BMM1}).
In \cite{HNR} we have very  recently revisited this problem with another approach, namely by developing dispersive $L^\infty$ linearized estimates in the physical space, which allowed us to use a  Lagrangian strategy in the spirit of \cite{BD} for the non-linear problem (see also \cite{Trinh}). In particular, \cite{HNR} shows that in the screened case, in all dimensions, the linear decay
 in the physical space is the same as for free transport, up to a logarithmic correction.

One expects the situation to be radically different for the unscreened Coulomb case~\eqref{linVP}, as evidenced in the pioneering works by Glassey and Schaeffer \cite{GS1,GS2}. In particular \cite{GS1,GS2} prove that in dimension $d=1$, when $\mu$ is a Maxwellian, the $L^2$ norm of the density of the solution to~\eqref{linVP} cannot in general decay faster than $1/ (\log t)^{13/2}$ (whereas for free transport it decays like $1/t^{1/2}$).
Furthermore, \cite{GS1,GS2} provide decay estimates, highlighting the influence of the rate of decay of $\mu$ at infinity:
\begin{itemize}

\item when $\mu$ is a Maxwellian, $\rho$ decays logarithmically fast in $L^2$ and $L^\infty$ norm.

\item When $\mu$ decays at most polynomially fast, $\rho$ decays polynomially fast in $L^2$ and $L^\infty$ norm (with a rate that cannot be better than $1/2$ and gets worse  when $\mu$ decays faster).

\item On the other hand when $\mu$ is compactly supported, they show that the $L^2$ norm of the density may not decay at all.

\end{itemize}

In this work, we shall consider a general class of  analytic  homogeneous equilibria (that includes Maxwellian and power laws for example). Quantitatively, we assume that there exist $R_{0}>0$ and $C_{0}>0$   so that for all polynomials $P$ of degree less than or equal to $\alpha_d$, with $\alpha_d := d+ 8$,
\begin{equation}
\label{muanal}
 |  \mathcal{F}_{v}( \mu)(\xi) | + | \mathcal{F}_{v}( P(v) \nabla_{v} \mu)(\xi) | \leq C_{0} e^{ - R_{0} | \xi |},
   \quad \forall \xi \in \mathbb{R}^d.
 \end{equation}
 where we use  $\mathcal{F}_{v}$ to denote the Fourier transform. %$\hat\cdot$ 

Following \cite{MV} and \cite{BMM2}, one could  expect that a relevant notion of stability is the one  in the sense of Penrose, that would correspond to asking that there is $\kappa>0$ such that
\begin{equation}
\label{Penrose-VP}
\inf_{\gamma \geq 0, \, \tau \in \R, \,  \xi \in \R^d}  \left| 1   - \int_{0}^{+ \infty} e^{-(\gamma + i \tau)s}\, {i \xi \over  |\xi|^2}\cdot  \mathcal{F}_v( \na_v \mu)(\xi s) \, d s \right| \geq \kappa.
\end{equation}
However, as we will soon observe,  though relevant in the torus case, this condition can \emph{never} be satisfied on the whole space. This is because of a \emph{low frequency} (in space) singularity (i.e. for small values of $|\xi|$), which is the reason why the decay that can be obtained in the screened case should not be expected here.  This  explains (most of) the results of \cite{GS1,GS2}: their strategy is based on a cut-off argument around the singularity, which accounts for why the rate of decay of $\mu$ matters in their result.
We shall show that despite this singularity, with a relevant notion of stability, a natural and much stronger  decay estimate that depends only on the 
dimension can be obtained.

A simplified version of our main result is stated in the following theorem.
\begin{thm}
\label{thm1}
%Assume that  \eqref{muanal}, \eqref{muradial} and (H1)-(H2) are satisfied.  
Let $d\geq 1$. Let $\mu$ be a non-negative radial equilibrium satisfying~\eqref{eq:norma} and ~\eqref{muanal}, of the form
$\mu(v) = F\left(\frac{|v|^2}{2}\right)$, with $F'(s) <0$, $\forall s \geq 0$.
%
% $\int_{\R^d} F'\left(\frac{|v|^2}{2}\right) \, dv \neq  0$. In addition, assume
%\begin{itemize}
%\item if $d=1$, $F' <0$, 
%\item if $d=2$, $\int_{\R^2} F'\left(\frac{|v|^2}{2}\right) \, dv <  0$,
%\item if $d=3$, $F>0$.
%\end{itemize}
Consider the density   $\rho(t,x)$ of  the solution of~\eqref{linVP}.
Then we can decompose
$$ \rho(t,x) = \rho^{R}(t,x) + \rho^S_{+}(t,x) + \rho^S_{-}(t,x),$$ 
where for all $t\geq 2$, we have
 $$   \| \rho^{R} (t)  \|_{L^\infty} \lesssim  \frac{\log t}{t^d} \left( \|f_{0}\|_{L^1_{x,v}} +   \|f_{0}\|_{L^1_{x}L^\infty_{v}}  \right)$$
  and for $k=0, \, 1$, 
 %  if   $d \geq 3$,  only  for  $k=1$ if  $d=1,2$,  we have 
$$   \|\rho^S_{\pm}(t) \|_{L^\infty} \lesssim \frac{\log t}{t^{{d \over 2} +k-1} } \sum_{0 \leq l \leq k} \left( 
\| \langle v \rangle^l \nabla_{v}^lf_{0}\|_{L^1_{x,v}} +   \| \langle v \rangle^l  \nabla_{v}^l f_{0}\|_{L^1_{x}L^\infty_{v}}  \right).
$$
%\begin{equation}
%\| \rho(t) \|_{L^\infty(\R^d)}  \lesssim \frac{1}{\langle t \rangle^{\frac{(d+1)^2}{4d+3}-1}}\left( \sup_{[0, t]} \,  \langle s\rangle^d\|S(s)\|_{L^\infty} + \sup_{[0, t]} \|S(s)\|_{L^1}\right),
%\end{equation}
%where $S(t,x)= \int_{\R^d} (\mathcal{S}(t,x,v) + f_0(x-tv,v)) \, dv$.

\end{thm}

\begin{rem}
All the assumptions,   in particular, $F'<0$,  are   satisfied   when $\mu$  is a  Maxwellian equilibrium  or  a power law
 $\mu(v)= c_{d} {1 \over (1 + |v|^2)^m}$, with $m$ sufficiently large.

\end{rem}

\begin{rem}
The rate of decay of $\mu$ does not play any role in this result.
\end{rem}

\begin{rem}
Observe that we do not state any decay in $L^2$ and therefore this is not in contradiction with \cite{GS1,GS2}.
\end{rem}

\begin{rem}
Contrary to the screened case, we do not have $\na_x^k \rho$ to decay faster than $\rho$; this is due to the singular part $\rho^S_\pm$.
\end{rem}

In Theorem~\ref{thm1} (and in  the more general version Theorem \ref{thm2} below), we have only stated $L^\infty$ type dispersive estimates. Nevertheless, we shall provide
 a much more precise description  of the structure of $\rho^{R}$ and $\rho^S_{\pm}$ in the following.
  By using  by now standard interpolation estimates, we could  deduce from them Strichartz estimates for example (\cite{Castella-Perthame}, \cite{Ginibre-Velo}).
  We observe that the regular part $\rho^{R}$ enjoys the same decay estimates as the solution of the linearized screened
  Vlasov-Poisson system obtained in \cite{HNR} which are themselves similar to the ones of the free transport up to the logarithmic factor. The singular part $\rho^S_{\pm}$ is precisely due to  a singularity at the frequency  
   $\xi= 0$, $\tau = \pm1$ in the dispersion relation. It can be seen as the solution to a dispersive  partial differential equation.
    Indeed, we shall show that $\rho^{S}_{\pm}$ is under the form
    $$   \rho^{S}_{\pm} (t,x)= \int_{0}^t G_{\pm}^S(t-s) *_{x} S(s)\, ds, \quad S(t)= \int_{\mathbb{R}^d} f_{0}(x-vt, v)\, dv$$
    the kernel $G_{\pm}^S$ being under the form
    $$  G^S_{+}(t,x) = \int_{\mathbb{R}^d} e^{Z_{\pm}(\xi) t + i x \cdot \xi} A_{\pm}(\xi) \, 
    d \xi$$
    where $ A_{\pm}$ is a smooth amplitude that is compactly supported for small $|\xi|$. The phase $Z_{\pm}(\xi)$ is such that $ \mbox{Re } Z_{\pm}(\xi) \leq 0$, 
    in addition,  $ \xi \mapsto \mbox{Re } Z_{\pm}(\xi)$ vanishes at $\xi=0$ and  is very flat  (and gets  flatter  when $\mu$ decays faster) so that only a very weak decay connected to the rate of decay of $\mu$  can be obtained
     from this piece of information. This accounts for the decay results of \cite{GS1,GS2}. %is in particular what is used in \cite{GS1}, \cite{GS2}. 
     Here we shall use that the imaginary part of the phase
      is non-degenerate so that  the decay rate $t^{- {d \over 2}}$ can be obtained from a stationary phase analysis.
       A significant part of the analysis of the paper  will be to perform a careful analysis of the  singularity of the dispersion relation at $\tau= \pm1$, $\xi= 0$
     and to justify that it  gives rise to the above singular term.

\section{Statement of the theorem with general assumptions on the equilibrium}

As a matter of fact, Theorem~\ref{thm1} is a special case of a more general result, allowing for a wider class of homogeneous equilibria
(not necessarily radial)
that satisfy a series of assumptions, which we now present.

\bigskip

\noindent{\sc Symmetry assumptions.} For all monomials $P$ of  odd degree $k \leq \alpha_d-1$,
we require that
\begin{equation}
\label{muodd}
\int_{\R^d} P(v) \mu(v) \, dv =0.
\end{equation}
We shall also ask that for all  $p \in \N\setminus\{0\}$ such that $2 p\leq  \alpha_d-1$, 
\begin{equation}
\label{muradial}
\exists C_{\mu}^p, \quad \forall \xi \in \R^d, \quad \int_{\mathbb{R}^d}  (\xi \cdot v)^{2p} \mu(v)\, dv  = C_\mu^p |\xi|^{2p}.
\end{equation}
Observe that~\eqref{muodd} is in particular satisfied when $\mu$ is even and~\eqref{muradial} when $\mu$ is radial; however both can also be satisfied assuming (many) algebraic identities on integrals of $\mu$ against polynomials.
 For most of the arguments, we shall  only need \eqref{muodd} and
 we will emphasize precisely where the additional assumption~\eqref{muradial} is needed in the paper.
 
  %We shall also   need a stronger version of symmetry for $\mu$:
% \begin{equation}
% \label{muradial}
% \mu \mbox{ is  radial}.
% \end{equation}

\bigskip

\noindent{\sc Stability assumptions.} Two stability assumptions are required.

\noindent {\bf Assumption (H1).} We shall  first ask for the stability condition: for every $\xi \neq 0$
\begin{equation}
\label{Penrose-VP-awayfrom0}
\inf_{\gamma \geq 0, \, \tau \in \R} \left| 1   -  \int_{0}^{+ \infty} e^{-(\gamma + i \tau)s}\, {i \xi \over  |\xi|^2}\cdot  \mathcal{F}_v( \na_v \mu)(\xi s) \, d s \right|>0,
\end{equation}
which is a weaker non-quantitative version of~\eqref{Penrose-VP} 
%a stronger version than  the Penrose stability condition on the torus, and is in particular satisfied for radial positive equilibria $\mu$. %satisfying  As a matter of fact, we will actually need a stronger version of this property (...)
% that is automatically satisfied as soon as $\mu $ is a radial positive function. We will therefore systematically make this latest assumption in this work.

\medskip

In order to tame the effect of the singularity at $\xi=0$, we shall require \emph{another Penrose stability} condition.
To this end, let us introduce 
$$ m_{KE}(z, \eta)=  -   \int_{0}^{+ \infty} e^{-(\gamma + i \tau)s}\, {i \eta \over  |\eta|^2}\cdot   \sum_{k,l} \eta_k \eta_l \mathcal{F}_v( v_k v_l \na_v \mu)(\eta s) \, d s, \quad z = \gamma + i \tau.$$
 For $\eta \in  \mathbb{S}^{d-1}$, thanks to \eqref{muanal}, we observe that $m_{KE}$ is holomophic in ${\operatorname{Re}}\, z >- R_{0}$.
 % We shall assume that: 
  
 \noindent {\bf Assumption (H2).}
  For every $\eta \in  \mathbb{S}^{d-1}$, there is only one zero of  $ z \mapsto 1 -  m_{KE}(z, \eta)$  on  $ {\operatorname{Re}}\, z =  0$
   which is $z= 0$. Moreover it verifies
  \begin{equation}
  \label{penrose3}
   \partial_{z} m_{KE}(0, \eta)= 0, \quad   \partial_{z}^2 m_{KE}(0, \eta)\neq 0,  \quad \forall \eta \in \mathbb{S}^{d-1}.
   \end{equation}  
This condition can be interpreted as a kind of Penrose stability condition for the so-called \emph{kinetic Euler} equation, which is a singular Vlasov equation arising in the quasineutral limit of the Vlasov-Poisson system and in Brenier's  incompressible optimal transport \cite{Br89}, \cite{BrCPAM}.

We are finally in position to state the main result of the paper.
\begin{thm}
\label{thm2}
Assume that  \eqref{eq:norma}, \eqref{muanal}, \eqref{muodd}, \eqref{muradial}, (H1) and (H2) are satisfied.
Then the conclusions of Theorem~\ref{thm1} hold.
\end{thm}

\begin{rem}The assumption (H2) can be replaced by

{\bf Assumption (H2').}  For every $\eta \in  \mathbb{S}^{d-1}$, there is no zero of  $ z \mapsto 1 -  m_{KE}(z, \eta)$  on  $ {\operatorname{Re}}\, z =  0$.

The proof of Theorem~\ref{thm2} gets slightly simplified in that case. However we have decided to focus on (H2) as (H2') is never satisfied for radial equilibria.
\end{rem}

Theorem~\ref{thm1} follows from Theorem~\ref{thm2} once that we have checked that the radial equilibria  $\mu= F(|v|^2/2)$ that we consider satisfy all required assumptions:
\begin{itemize}

\item we have already seen that~\eqref{muodd} and~\eqref{muradial} are satisfied when $\mu$ is radial.

\item 
As seen from \cite[Proposition 2.1 and Remark 2.2]{MV}, the assumption (H1) is verified for radial equilibria in any dimension assuming
that  $F' <0$.

\item  Finally the assumption (H2) is 	also  satisfied if $F'<0$. We postpone the proof of this result to an appendix, see Section~\ref{appendix}.

\end{itemize}

The rest of the paper is dedicated to the proof of Theorem~\ref{thm2}. %We denote by $\mathcal{F}$ the Fourier transform in time and space, and by $\mathcal{F}_v$ the Fourier transform in variable $v$. 

\section{Reduction to kernel estimates}
We study the linear equation
\begin{equation}\label{def-m}
\rho(t,x) = \int_0^t \int_{\R^d}  -  [\na_x  \Delta_x^{-1} \rho](s,x -(t-s)v) \cdot \na_v \mu(v) \, dv ds + S(t,x), \quad t \geq 0,
\end{equation}
with $S$ being a given source term (as we shall see later,  we can rewrite \eqref{linVP} under this form by integrating along
 the characteristics of the free transport).  In what follows, we extend $\rho$ and $S$ by zero for $t<0 $ so that the equation \eqref{def-m}
is satisfied for $t \in \mathbb{R}$. 
For $\gamma >0$ sufficiently large, by using the Fourier transform in space and time, we get that the solution of \eqref{def-m}
is given by 
$$  \mathcal{F}(e^{-\gamma t } \rho) (\tau, \xi) =  m_{VP}(\gamma, \tau, \xi)   \mathcal{F}(e^{-\gamma t} S)(\tau, \xi),$$
where $\mathcal{F}$ denotes the Fourier transform in time and space,
 and hence that
 \begin{equation}
 \label{eqrho} \rho(t,x) =  S +  \left( e^{\gamma  t }\left( \mathcal{F}^{-1} {  m_{VP}(\gamma, \cdot) \over 1 - m_{VP}(\gamma, \cdot)} \right)\right)*_{t, x} S = S+ G*_{t,x}S,
 \end{equation}
where 
 \begin{equation}
 \label{Gdef}G(t,x) = \int_{\mathbb{R} \times \mathbb{R}^d} e^{ \gamma t + i \tau t} e^{i x \cdot \xi}
   {  m_{VP}(\gamma, \tau, \xi) \over 1 - m_{VP}(\gamma, \tau, \xi)}\, d\tau d\xi.
   \end{equation}
The aim of the remaining will be to estimate the kernel $G$. Note that the definition of the kernel  $G$ depends on $\gamma$, but that in regions
where the integrand is an holomorphic function of $z= \gamma + i \tau$ it actually does not depend on $\gamma$ since  we can appropriately change the integration contour, via the Cauchy formula, without
 changing $G$. In particular, by taking the limit $\gamma \rightarrow + \infty$, we get that $G_{|t<0}=0.$
 
% For the use of this argument in the torus case, we refer to \cite{MV} (see also \cite{Degond}, \cite{GNR}).
%  where pointwise bounds on the resolvent and an exponential decay in time for the kernel were derived. 

 % For the use of this argument in the torus case, where an exponential decay in time can be obtained, 
% we refer to \cite{MV} (see also \cite{Degond}, \cite{GNR}).

Precisely, in this paper, we shall prove: 
 \begin{thm}
 \label{theokernel}
 For $t \leq 1$, we have the estimate
 $$ \|G(t) \|_{L^\infty} \lesssim { 1 \over t^{d-1}}, \quad \|G(t) \|_{L^1} \lesssim t.$$
 Moreover, 
 assuming \eqref{eq:norma}, \eqref{muanal}, \eqref{muradial}, (H1), (H2), we can write for $t \geq 1$,
 $$ G= G^{R}(t,x) + G^S_{+}(t,x)+ G^S_{-}(t,x), $$
 where 
 $$  \|G^{R}(t) \|_{L^\infty} \lesssim { 1 \over t^{d+1}}, \quad \|G^{R}(t) \|_{L^1} \lesssim {1 \over t}, \quad \forall t \geq 1$$
 and 
 $$   \|G^{S}_{\pm}(t) \|_{L^\infty} \lesssim { 1 \over t^{d\over 2}}, \quad \|G^{S}_{\pm}(t) \|_{L^2} \lesssim 1, \quad \forall t \geq 1.$$
  \end{thm}
A more accurate description of $G^S_\pm$ is given in Proposition \ref{propGS}. They  can be seen as the kernel of the propagator
 of a dispersive PDE.
\section{Properties of the symbol $m_{VP}$}

Let us set for $(\gamma, \tau, \xi) \in \R \times \R \times \R^d$,
\begin{align}
 m_{VP}(\gamma, \tau, \xi)&= \int_{0}^{+ \infty} e^{-(\gamma + i \tau)s}\, {i \xi \over  |\xi|^2}\cdot  \mathcal{F}_v( \na_v \mu)(\xi s) \, d s, \\
  m_{VB}(\gamma, \tau, \xi)&= \int_{0}^{+ \infty} e^{-(\gamma + i \tau)s}\, {i \xi}\cdot  \mathcal{F}_v( \na_v \mu)(\xi s) \, d s, \\
  m_{KE}(\gamma, \tau, \xi)&= - \int_{0}^{+ \infty} e^{-(\gamma + i \tau)s}\, {i \xi \over  |\xi|^2}\cdot   \sum_{k,l} \xi_k \xi_l \mathcal{F}_v( v_k v_l \na_v \mu)(\xi s) \, d s.
\end{align}

\begin{rem}
As already mentioned, $m_{KE}$ is the symbol associated to the Kinetic Euler equation.
It turns out that the symbol 
 $m_{VB}$ is the one associated to the so-called \emph{Vlasov-Benney} equation (see e.g. \cite{Bardos}) which is another singular Vlasov equation that shows up in the quasineutral limit of the Vlasov-Poisson system \cite{HKR}.
\end{rem}

 Let us  define $\Omega_{R_{0}} = \mathcal{A} \cap \mathcal{C}_{R_{0}}$ where, 
  $$\mathcal{A}=\left\{ (\gamma, \tau,  \xi) \in \mathbb{R}^{d+2}, \, { 1 \over 2} < | \gamma |+ | \tau | + |\xi| < 2 \right\}, 
   \quad   \mathcal{C}_{R_{0}}= \left\{ (\gamma, \tau , \xi) \in \mathbb{R}^{d+2}, \, \xi \neq 0, \,  \gamma >- R_{0} | \xi|\right\}.  $$
 
 \subsection{Estimates of $m_{VP}$}
The following is an adaptation of Lemma 2.2 in \cite{HNR}. We get stronger properties due to the regularity assumption
 \eqref{muanal} and the symmetry assumption~\eqref{muodd}.
\begin{prop}
  \label{lem-VP-VB}
  Assuming \eqref{eq:norma}, \eqref{muanal} and \eqref{muodd}, we have the following properties.
For every  $(\gamma, \tau, \xi) \in \mathcal{C}_{R_{0}}$,
      \begin{equation}
     \label{eq-K1def} m_{VP}(\gamma, \tau, \xi)=    \frac{1}{|\xi |^2}   m_{VB}(\gamma, \tau, \xi), \quad 
      m_{VP}(\gamma, \tau, \xi)= - { 1 \over (\gamma + i \tau)^2}( 1 - m_{KE}(\gamma, \tau, \xi))
     \end{equation}
The symbols $  m_{VB}(\gamma, \tau, \xi)$  and $m_{KE}(\gamma, \tau, \xi)$ 
 are for $\xi \neq 0$ holomorphic with respect to the variable $z= \gamma + i \tau$ in $\gamma > - R_{0} | \xi|$.
  Moreover, they 
are  positively homogeneous of degree zero and $m_{VB},\, m_{KE} \in \mathscr{C}^{\alpha_d - 2}(\mathcal{C}_{R_{0}/2})$.
 Quantitatively, there exists $C>0$ such that 
   \begin{equation}
   \label{eq-estmVB}
   | \partial_{z}^\alpha \partial_{\xi}^\beta m_{VB}(\gamma, \tau, \xi) | + | \partial_{z}^\alpha \partial_{\xi}^\beta m_{KE}(\gamma, \tau, \xi) |
    \leq {C \over |(\gamma, \tau, \xi)|^{|\alpha| + |\beta|}}, \, \forall | \alpha| + |\beta | \leq   \alpha_d - 2,\,  \, \forall (\gamma, \tau, \xi) \in  \mathcal{C}_{R_{0}/2},
    \end{equation}
where we use $\partial_{z}=  \partial_{\gamma} - i \partial_{\tau}.$
   
  \end{prop}
 
\begin{rem}
In the following, we shall often abuse notations and write the symbols as functions of $(\gamma, \tau,\xi)$ or $(z,\xi)$ depending on what is most convenient.
\end{rem}

%We continue by studying the symbol $m_{KE}$ on the sphere.
\begin{proof} 
Let us first observe that thanks to \eqref{muanal}, $m_{VP}$, $m_{KE}$ and $m_{VB}$ are well-defined
 in $\mathcal{C}_{R_{0}}$ and holomorphic in $z$ for $\xi\neq 0$ and $\operatorname{Re}\, z > - R_{0} | \xi|$.
 Let us prove  \eqref{eq-K1def}. The first relation is trivial. For the second one, 
by two successive integrations by parts in $s$, we obtain
\begin{multline*}
m_{VP}(\gamma, \tau, \xi)  = \frac{1}{\gamma +i \tau}  \int_{0}^{+ \infty} e^{-(\gamma + i \tau)s}\, {i \xi \over  |\xi|^2}\cdot  \pa_s \left(\mathcal{F}_v( \na_v \mu)(\xi s)\right) \, d s  %+ \frac{1}{\gamma +i \tau} \left[ \right]_0^{+\infty} 
\\ =   \frac{1}{(\gamma +i \tau)^2}  \int_{0}^{+ \infty} e^{-(\gamma + i \tau)s}\, {i \xi \over  |\xi|^2}\cdot  \pa^2_{s} \left(\mathcal{F}_v( \na_v \mu)(\xi s)\right) \, d s  \\
  + \frac{1}{(\gamma +i \tau)^2}  {i \xi \over  |\xi|^2} \cdot \pa_s \left(\mathcal{F}_v( \na_v \mu)(\xi s)\right)|_{s=0}.
\end{multline*}
Since by~\eqref{eq:norma}, 
$$
 \pa^2_{s} \left(\mathcal{F}_v( \na_v \mu)(\xi s)\right) = -  \sum_{k,l} \xi_k \xi_l \mathcal{F}_v( v_k v_l \na_v \mu)(\xi s) , \quad
 {i \xi \over  |\xi|^2} \cdot  \pa_s \left(\mathcal{F}_v( \na_v \mu)(\xi s)\right)|_{s=0} = -1,
$$
we finally get   \eqref{eq-K1def}. The degree zero homogeneity property comes from a  straightforward change of variable.
It  remains to prove \eqref{eq-estmVB}. We first give the proof  for $m_{VB}$.

 Since  we have 
  $$  m_{VB}(\gamma, \tau, \xi)=  \int_{0}^{+\infty} e^{- (\gamma+ i \tau )t}   i  \xi \cdot \widehat{\nabla_{v} \mu}
  (t \xi) \, dt, $$
   we get by using  \eqref{muanal}, that 
         $$  |m_{VB}(\gamma, \tau, \xi)| \leq C \int_{0}^{+ \infty} | \xi | e^{-  R_{0}  | \xi|  t/2 } \, dt \leq C, 
          \quad \forall (\gamma, \tau, \xi) \in \mathcal{C}_{R_{0}/2} .$$
  We now estimate the derivatives in $ \Omega_{R_{0}/2}$;  let us first handle the case when $|\xi| \geq {1 \over 4}.$ 
Thanks to \eqref{muanal},  we also  have that 
   $$ | \partial_{z}^\alpha \partial_{\xi}^\beta m_{VB}(\gamma, \tau, \xi) | \lesssim 
     \int_{0}^{+ \infty}    \langle t \rangle^{|\alpha| + | \beta|}  | \xi | e^{-  R_{0}  | \xi|  t/2 } \, dt $$
     and therefore, for $|\xi| \geq {1\over 4}$ and $| \alpha| + |\beta | \leq \alpha_d$, we obtain
     \begin{equation}\label{bd-X1}   | \partial_{z}^\alpha \partial_{\xi}^\beta m_{VB}(\gamma, \tau, \xi) | \lesssim  1, \qquad  (\gamma, \tau, \xi) \in \Omega_{R_{0}/2}, \, 
\quad      | \xi | \geq { 1 \over 4}.\end{equation}

Let us next consider the case  $|\xi| \leq { 1 \over 4}$, in which we make use of the fact that $|z|$ is positively  bounded from below, recalling $(\gamma,  \tau, \xi) \in \mathcal{A}$. Integrating by parts again, we get for every $n = 2,\cdots, \alpha_d$,
   \begin{equation}
   \label{Kxismall}
    m_{VB}(\gamma, \tau, \xi)=  \sum_{k= 2}^n
     { 1 \over z^k}  \mathscr{P}_{k}(\xi)   +   { 1 \over z^n} R_{n}(\gamma, \tau, \xi)
    \end{equation} 
%     { 1 \over (i\tau)^n} \left( \int_{0}^{+ \infty} e^{-i \tau t} i \xi\cdot (D_{\xi}^k \widehat{\nabla_{v} \mu}) (t\xi)\, dt \right):\xi^{\otimes k}
%     \end{equation}
     where
     $$\begin{aligned} \mathscr{P}_{k}(\xi)&= (-1)^{k-1}i^k \xi \cdot \mathcal{F}\left( v^{\otimes k-1}\nabla_{v} \mu \right) (0): \xi^{\otimes k-1}, \\
       R_{n}(\gamma, \tau, \xi) &=  \int_{0}^{+ \infty} e^{-(\gamma + i \tau) t} r_{n}(t,\xi)\, dt, \quad r_{n}(t,\xi)=
      (-1)^n i^n \xi\cdot \mathcal{F}\left(v^{\otimes n}\nabla_{v} \mu\right) (t\xi):\xi^{\otimes n},
       \end{aligned}
       $$
        with the definition
  $$ \xi \cdot  \mathcal{F}(v^{\otimes k} \nabla_{v} \mu)(\zeta): \xi^{\otimes k}= \sum_{j_{0}, j_{1}, \cdots j_{k}}
  \xi_{j_{0}}  \xi_{j_{1}} \cdots \xi_{j_{k}} \mathcal{F}( v_{j_{1}}\cdots v_{j_{k}} \partial_{v_{j_{0}}} \mu)(\zeta).$$
     Note that $\mathscr{P}_{k}$ is a homogeneous polynomial of degree $k$.  Thanks to \eqref{muanal}, we have
   $$  | r_{n}(t, \xi) |
    \lesssim  |\xi|^{n+1} e^{ - R_{0} | \xi| t }.$$
 More generally,  we have for all $|\beta|\le n$, 
    $$  |\partial_{\xi}^\beta  r_{n}(t, \xi) | \leq  |\xi|^{n+1- |\beta |} e^{ - R_{0} | \xi| t }.$$

   Consequently, applying derivatives  to the expansion  \eqref{Kxismall} and using the above estimates with $n=\alpha_d$, we get for $|\xi| \leq { 1 \over 4}$
    and $(\gamma, \tau, \xi) \in \Omega_{R_{0}/2}$ (which in particular implies that $|z|\ge \frac14$), 
    $$  | \partial_{z}^\alpha \partial_{\xi}^\beta m_{VB}(\gamma, \tau, \xi) | \lesssim 
     1 +  \int_{0}^{+ \infty}  t^{|\alpha|} |\xi|^{\alpha_d+1- | \beta|} e^{ - R_{0} | \xi| t /2}\, dt
      \lesssim 1 +  \int_{0}^{+ \infty}  s^{|\alpha|} |\xi|^{\alpha_d + 1- | \beta| - |\alpha |}  e^{ - R_{0} s /2}\, ds.$$
 Thus, we get for $| \alpha| + |\beta | \leq   \alpha_d$,
   $$  | \partial_{z}^\alpha \partial_{\xi}^\beta m_{VB}(\gamma, \tau, \xi) | \lesssim 
     1, \qquad  (\gamma, \tau, \xi) \in \Omega_{R_{0}/2}, \, \quad
      | \xi | \leq { 1 \over 4}.$$        
      This, together with \eqref{bd-X1}, concludes the proof of the estimates for $m_{VB}$ on $\Omega_{R_0/2}$;
      we finally obtain~\eqref{eq-estmVB} for $m_{VB}$ by degree zero homogeneity.
      
    Let us now prove the estimates for $m_{KE}$ on $\Omega_{R_0/2}$. The same argument as above applies for $| \xi| \geq {1/4}$.
   Thus it suffices to study $| \xi| \leq {1/4}$. As before, we can integrate by parts to get that
         \begin{equation}
   \label{mKExismall}
    m_{KE}(z, \xi)=  \sum_{k= 2}^n
     { 1 \over z^k}  \mathcal{Q}_{k}(\xi)   +   { 1 \over z^n} R^{KE}_{n}(\gamma, \tau, \xi)
    \end{equation}
    where
    $$
  \begin{aligned} \mathcal{Q}_{k}(\xi)&= (-1)^{k}i^k {\xi \over |\xi|^2} \cdot \mathcal{F}\left( v^{\otimes k+1}\nabla_{v} \mu \right) (0): \xi^{\otimes k+1}, \\
       R_{n}^{KE}(\gamma, \tau, \xi) &=  \int_{0}^{+ \infty} e^{-(\gamma + i \tau) t} r_{n}^{KE}(t,\xi)\, dt, \quad r_{n}^{KE}(t,\xi)=
      (-1)^{n+1} i^n {\xi \over | \xi|^2}\cdot \mathcal{F}\left(v^{\otimes n+2}\nabla_{v} \mu\right) (t\xi):\xi^{\otimes n+2}.
       \end{aligned}
 $$  
 In this case, we need to study more carefully the structure of this expansion since  the function $\xi \otimes \xi/ | \xi|^2$
  multiplied by a polynomial of $\xi$
  is not necessarily a smooth function of $\xi$.  We first observe that if $k$ is odd, we have
  $$  \mathcal{Q}_{k}(\xi) = 0, \quad \forall \xi.$$
  Indeed,  we  can integrate by parts and use that 
  $$ \int_{\mathbb{R}^d}  v^{\otimes k} \mu(v) \, dv = 0$$
   if $k$ is odd thanks to the symmetry assumption~\eqref{muodd}.  We thus have the expansion
    \begin{equation}
   \label{mKExismallbis}
    m_{KE}(z, \xi)=  \sum_{p= 1}^l
     { 1 \over z^{2p}}  \mathcal{Q}_{2p}(\xi)   +   { 1 \over z^{2l+1}} R^{KE}_{2l+1}(\gamma, \tau, \xi).
    \end{equation}  
Moreover,  we have
\begin{multline*} 
\mathcal{Q}_{2p}(\xi) =  (-1)^p {\xi \over | \xi |^2} \cdot  \mathcal{F}\left( v^{\otimes 2p+1}\nabla_{v} \mu \right) (0): \xi^{\otimes 2p+1}
\\=  (-1)^p { 1 \over | \xi|^2} \sum_{j_{0}, \cdots, j_{2p+1}} \xi_{j_{0}} \cdots \xi_{j_{2p+1}} \int_{\mathbb{R}^d}
 v_{j_{1}}\cdots v_{j_{2p+1}} \partial_{v_{j_{0}}} \mu(v) \, dv.
 \end{multline*}
 Integrating by parts, we observe that if $j_{1}, \cdots,  j_{2p+1}$ are all different from $j_{0}$, the integral vanishes, 
  therefore, after relabelling, we get that
  \begin{multline*} 
\mathcal{Q}_{2p}(\xi) =  
   (-1)^{p+1} (2p+1) { 1 \over | \xi|^2} \sum_{j_{0}, j_{2}\cdots, j_{2p+1}} (\xi_{j_{0}})^2 \xi_{j_{2}} \cdots \xi_{j_{2p+1}} \int_{\mathbb{R}^d}
 v_{j_{2}}\cdots v_{j_{2p+1}} \mu(v) \, dv \\
 = (-1)^{p+1} (2p+1)\sum_{ j_{2}, \cdots,  j_{2p+1}}\xi_{j_{2}} \cdots \xi_{j_{2p+1}} \int_{\mathbb{R}^d}
 v_{j_{2}}\cdots v_{j_{2p+1}} \mu(v) \, dv 
 \\=  (-1)^{p+1} (2p+1)  \int_{\mathbb{R}^d}( \xi \cdot v)^{2p} \mu(v)\, dv.
 \end{multline*}
 We have thus obtained that $\mathcal{Q}_{2p}$ is a polynomial in $\xi$ and hence a smooth function of $\xi$.  
 We can then get from   the expansion \eqref{mKExismallbis}
  and similar estimates as for $m_{VB}$  that the derivatives with respect to  $z$ and $\xi$ of order less than $\alpha_d - 2$
  are uniformly bounded for $| \xi| \leq 1/4$. Finally, the estimate ~\eqref{eq-estmVB} follows by using the degree zero homogeneity.
    \end{proof}

    In the proof of Proposition~\ref{lem-VP-VB}, we have obtained a refined asymptotic expansion of $m_{KE}$, for which we have relied on the symmetry assumption~\eqref{muodd}. We gather this very useful statement in the following lemma.
    %also needed to handle a singularity near $\xi =0$.
    
  \begin{lem}
  \label{remarkmKE}
%  Note that we have established the following useful refined expansion of $m_{KE}$
  Assuming \eqref{eq:norma}, \eqref{muanal} and \eqref{muodd}, we have the following  expansion of $m_{KE}$ for all $l \leq \lfloor \frac{\alpha_d-3}{2} \rfloor$:  
  \begin{align}
  \label{mKErefined}
  & m_{KE}(z, \xi)=  \sum_{p= 1}^l
     { 1 \over z^{2p}}  \mathcal{Q}_{2p}(\xi)   +   { 1 \over z^{2l+1}} R^{KE}_{2l+1}(\gamma, \tau, \xi), \\
 & \label{Q2p}  \mathcal{Q}_{2p}(\xi) =  (-1)^{p+1} (2p+1) \xi^{\otimes 2p}  : \int_{\mathbb{R}^d} v^{\otimes 2p} \mu(v)\, dv
 = (-1)^{p+1} (2p+1)   \int_{\mathbb{R}^d}  (\xi \cdot v)^{2p} \mu(v)\, dv ,
  \\
&    \label{resteKE}  R_{2l+1}^{KE}(z, \xi) =  \int_{0}^{+ \infty} e^{-(\gamma + i \tau) t} r_{2l+1}^{KE}(t,\xi)\, dt, \,
  r_{2l+1}^{KE}(t,\xi)=
       (-1)^{l+1}    {\xi \over | \xi|^2}\cdot \mathcal{F}\left(v^{\otimes 2l+2}\nabla_{v} \mu\right) (t\xi):\xi^{\otimes 2l+3}
  \end{align}
  where the remainder  satisfies uniformly  for $(\gamma, \tau, \xi) \in \mathcal{C}_{R_{0/2}}$ the estimate
%{\color{red} $$XXX:.....R_{2l+1}^{KE}.... ?$$} 
 \begin{equation}
  \label{remindermKE}  |\partial_{z}^\alpha \partial_{\xi}^\beta R_{2l+1}^{KE}(z, \xi) |
   \lesssim  | \xi|^{2 l+ 1 - | \alpha | - | \beta |}, \quad | \alpha | + | \beta | \leq 2l +1.
  \end{equation}
  In particular, we get
  $$ \mathcal{Q}_{2}(\xi)=  - 3 H_{\mu} \xi \cdot \xi, \quad H_{\mu}= \int_{\mathbb{R}^d} v\otimes v\, \mu(v) \, dv. $$
  Under the additional symmetry assumption \eqref{muradial}, we have  
 \begin{equation}
 \label{Q2pbetter}   \mathcal{Q}_{2p}(\xi)= (-1)^p (2p+ 1 ) | \xi|^{2p} C^p_{\mu}
\end{equation}
 where
 $$  C^p_{\mu} = \int_{\mathbb{R}^d} v_{1}^{2p} \mu(v)\, dv.$$
 In particular, we have  
 $$ C_{\mu} := C^1_{\mu}= { 1 \over d} \int_{\mathbb{R}^d} |v|^2 \, \mu(v)\, dv.$$
   
  \end{lem}

As a consequence of Proposition~\ref{lem-VP-VB}, we obtain estimates for $m_{VP}$. 
   \begin{coro}
   \label{lem-symbolVP-sphere}
   The symbol $m_{VP}(\gamma, \tau, \xi)$ 
 for $\xi \neq 0$ is holomorphic with respect to the variable $z= \gamma + i \tau$ in $\operatorname{Re} z > - R_{0} | \xi|$.
  Moreover,  it 
is  positively homogeneous of degree  $-2$ and $m_{VP} \in \mathscr{C}^{\alpha_d - 2}(\mathcal{C}_{R_{0}/2})$.
 Quantitatively, there exists $C>0$ such that 
   \begin{equation}
   \label{eq-estmVP}| \partial_{z}^\alpha \partial_{\xi}^\beta m_{VP}(\gamma, \tau, \xi) | 
    \leq {C \over |(\gamma, \tau, \xi) |^{2 +|\alpha |+ |\beta|}}, \, \forall   | \alpha| + |\beta | \leq  \alpha_d - 2,\,  \, \forall (\gamma, \tau, \xi) \in \mathcal{C}_{R_{0}/2}.
    \end{equation}

    \end{coro}

\begin{proof}
The fact that $m_{VP}$ is positively homogeneous of degree $-2$ follows from a change of variables.
For the other properties, by using the previous lemma, it suffices to observe that
  \begin{equation}
  \label{mVPbig}   m_{VP}(\gamma, \tau, \xi)=  { 1  \over 8   |\xi|^2 - ( \gamma + i \tau )^2} \left( 8 m_{VB}(\gamma, \tau, \xi) + ( 1 - m_{KE}(\gamma, \tau, \xi) \right).
  \end{equation}
   We get  that  the modulus of the denominator  for $(\gamma, \tau, \xi)
      \in \Omega_{R_{0}/2}$ is positively  uniformly bounded from below. Indeed, for $|\xi|^2\geq 1/4$, we observe that 
$$
     \left| 8 |\xi|^2 - ( \gamma + i \tau )^2  \right| \geq 8 | \xi|^2 - \gamma^2
    >c_0,  $$
    for some $c_0>0$.
    Otherwise, if $|\xi|^2< 1/4$, we must have $\tau^2 + \gamma^2> 1/4$, in which case 
    $$
     \left| 8 |\xi|^2 - ( \gamma + i \tau )^2  \right| \geq |\tau^2- \gamma^2|  + 2 |\tau| |\gamma| 
    >c_1,   $$
    for some $c_1>0$.
%       \begin{equation}
%  \label{mVPbig}   m_{VP}(\gamma, \tau, \xi)=  { 1  \over R_{0}^2   |\xi|^2 - ( \gamma + i \tau )^2} \left( R_{0}^2 m_{VB}(\gamma, \tau, \xi) + ( 1 - m_{KE}(\gamma, \tau, \xi) \right).
%  \end{equation}
%   We get  that  the modulus of the denominator  for $(\gamma, \tau, \xi)
%      \in \Omega_{R_{0}/2}$ is bounded from below. Indeed, we have that 
%$$
%     \left| R_{0}^2   |\xi|^2 - ( \gamma + i \tau )^2  \right| \geq { 3 \over 4} R_{0}^2 | \xi|^2 + \tau^2 
%    \geq { 1  \over 2} R_{0}^2 | \xi |^2 + \tau^2 + \gamma^2 
%     \geq \min({ R_{0} \over 2}, 1) |(\gamma, \tau, \xi) |^2  
%     $$
   Consequently, we can apply~\eqref{eq-estmVB} and the estimates~\eqref{eq-estmVP} follow since $m_{VP}$ is  positively homogeneous of degree $-2$.
\end{proof}

\subsection{Zeroes of $1-m_{VP}$}

In this section we give a sharp description of the zeroes of $1-m_{VP}$.
As we shall see, they are localized in the region $|\gamma | \leq \eps_{3} | \xi|$, $ | \tau \pm 1| \leq \eps_3|\xi| $ and $|\xi | \leq \eps_{3}$ for some small $\eps_3$. Using the implicit function theorem, we are able to describe them by smooth curves. 

\begin{prop}
\label{lem-zeroes}
 Assuming that  \eqref{eq:norma}, \eqref{muanal} and \eqref{muodd} hold, we have the following properties:
 \begin{enumerate}
 \item[i)] There exists $M>0$ such that for every $(\gamma, \tau, \xi) \in \mathcal{C}_{R_{0}/2}$
  and $|(\gamma, \tau, \xi)| \geq M$, we have
  $$ |1-m_{VP}(\gamma, \tau, \xi) |\geq { 1 \over 2}.$$
  \item[ii)] Assuming  (H1), for every $\delta >0$, there exists $c_{\delta}>0$ and $R_{\delta}\in (0, R_{0}/2]$ such that
   for every $(\gamma, \tau, \xi) \in \mathcal{C}_{R_\delta}$ with $| \xi | \geq \delta$, we have 
  $$|1 - m_{VP}(\gamma, \tau, \xi)| \geq  c_{\delta}.$$
 \item[iii)] Assuming  (H2), there exists $R_{1} \in (0, R_{0}/2]$,  $\eps_{1}>0$  and $C_{\eps_{1}}>0$, $c_{\eps_{1}}>0$  such that for every 
 $(\gamma, \tau, \xi) \in \mathcal{C}_{R_{1}}$ and   $|\gamma| \leq  \eps_{1} |\xi|$, $ |(\tau, \xi)| \leq \eps_{1}$, we have
 \begin{equation}
 \label{penrosepetit}\left|  { m_{VP}(\gamma, \tau, \xi) \over 1 - m_{VP}(\gamma, \tau, \xi)} \right|  \leq  C_{\eps_{1}}, 
  \quad  |z^2 + 1 - m_{KE}(z, \xi) |  \geq c_{\eps_{1}} \min(1, |\tilde z|^2), \quad z= |\xi| \tilde z.
  \end{equation}
\item[iv)] There exists $\eps_{2}\in (0,  \min(R_{0}/2, \eps_{1})]$, $A_{0}\geq 1$,  such that for every
 $ \eps \in (0, \eps_{2}]$, for every $A \geq A_{0}$ and  for every  $(\gamma, \tau, \xi) \in \mathcal{C}_{R_1}$, with $| \xi| \leq \eps, $
 $ | \gamma | \leq  \eps | \xi| $,  and  $  | |\tau|  -   1 | \geq  A \eps | \xi|$,  ${1 \over \eps_{1}} \geq |\tau| \geq \eps_{1}$,
 %(where $\eps_{1}$ is given by iii)),  
 there holds
 \begin{equation}
 \label{belowmiddle}|1 - m_{VP}(\gamma, \tau, \xi) | \geq  A\,  \eps\, \eps_{1}^2 | \xi|/4.
 \end{equation}
 \item[v)]  Assuming (H1), there exists $\eps_{3}>0$, $\eps_{3}\in (0,  \min(R_{0}/2, \eps_{2})]$ such that for every $\xi \neq 0$, the zeroes of $1-m_{VP}(\gamma, \tau, \xi)$
  with  $|\xi | \leq \eps_{3}$, $|\gamma | \leq \eps_{3} | \xi|$, $ | |\tau| - 1 | \leq \eps_3|\xi| $ are given by two $\mathscr{C}^1$  curves
     $$ Z_{\pm}(r, \omega) = \pm i + r  \Gamma_{\pm}(r, \omega)    +  i r \mathcal{T}_{\pm}(r, \omega)$$ where
     $\xi= r \omega$, $\omega \in \mathbb{S}^{d-1}$, 
      $\Gamma_{\pm} \leq 0$,  $\Gamma_{\pm}(0, \omega)= 0 $,  $\partial_{r} \Gamma_{\pm}(0, \omega)=0,$ $\Gamma_{\pm}(r, \omega)<0$ for $r \neq 0$ and 
       $\mathcal{T}_{\pm}$ is real with  $ \mathcal{T}_{\pm}(0, \omega)=0$,     $\partial_{r} \mathcal{T}_{\pm}(0, \omega)\neq 0$.
 
 \end{enumerate}

\end{prop}

\begin{proof}
Let us start with i).
 To this end, we can apply Corollary \ref{lem-symbolVP-sphere}.  This entails that 
 $$| m_{VP} (\gamma, \tau, \xi) | \leq { C \over |(\gamma, \tau, \xi) |^2}, \quad \forall (\gamma, \tau, \xi) \in \mathcal{C}_{R_{0}/2}$$
and hence 
$$|1-m_{VP}(\gamma,\tau,\xi) | \geq { 1 \over 2}$$
if $|(\gamma, \tau, \xi)|$ is sufficiently large.

Let us prove ii).  By using i), the estimate is true if we have in addition $|(\gamma, \tau, \xi)| \geq M$, it thus suffices
 to consider the case that $ |\xi| \geq  \delta$ and  $|(\gamma, \tau, \xi)| \leq M$.
  By (H1) and by compactness, we get that
   $$ | 1 - m_{VP}(\gamma, \tau, \xi)| \geq  2 c_{\delta}$$
   for some $c_{\delta}>0$ if $\gamma \geq 0$. By continuity, the inequality without the factor $2$ remains true for $\gamma  = -  \alpha | \xi|$ for $\alpha \leq  c$
    with $c>0$ small enough.

To prove iii), we observe that 
by Proposition~\ref{lem-VP-VB}, we can write
\begin{equation}
\label{expiii2} {m_{VP} \over 1- m_{VP}} = - { 1 \over z^2} { 1- m_{KE} \over  1 + { 1 \over z^2}(1-m_{KE})}
 =  { 1 -m_{KE} \over  z^2 + 1- m_{KE}}.
 \end{equation}
   By degree zero homogeneity of $m_{KE}$, we can set $\tilde z =z/ |\xi|$, $\eta= \xi/|\xi|$ , with $\tilde z= \tilde \gamma + i \tilde \tau$
    and  $| \tilde \gamma| \leq \eps_{1}$, $|\tilde \tau| \leq \eps_{1}/|\xi|$.
  This yields
$${m_{VP} \over 1- m_{VP}} (z, \xi)=  { 1- m_{KE}(\tilde z, \eta) \over  |\xi|^2 \tilde z^2  +  1- m_{KE}(\tilde z, \eta)}.$$ 
 By using  \eqref{mKErefined}, we have that  uniformly for  $| \tilde \gamma| \leq  R_{0}/2$ and $\eta \in \mathbb{S}^{d-1}$, 
 $$ \lim_{|\tilde \tau| \rightarrow + \infty} |m_{KE}(\tilde z, \eta)| = 0.$$
   Therefore, for $|\tilde \tau| \geq M$ sufficiently large 
   $$ |1- m_{KE}(\tilde z, \eta)| \geq {1 \over 2}$$
   and hence for $\eps_{1}$  sufficiently small, we get
$$\left|  |\xi|^2 \tilde z^2  +  1- m_{KE}(\tilde z, \eta)  \right| = \left|  z^2  +  1- m_{KE}(\tilde z, \eta)  \right|  \geq {1 \over 4}.$$
As a consequence, we conclude that  $\left|{m_{VP} \over 1- m_{VP}}\right|$ is bounded.

In a similar way, for every $\tilde \eps >0$ if $ \tilde \eps \leq |\tilde  \tau | \leq M$, $1-m_{KE}(i\tilde \tau, \eta)$ does not vanish
thanks to (H2). By compactness and continuity this remains  true uniformly for $\tilde \gamma$ sufficiently small
 and $\eta \in \mathbb{S}^{d-1}.$ In this regime, we thus also get 
  $\left|  |\xi|^2 \tilde z^2  +  1- m_{KE}(\tilde z, \eta)  \right|$ is uniformly strictly positive  
 and  also   that  ${m_{VP} \over 1- m_{VP}}$ is bounded. 
   
   Consequently there only remains to study the vicinity of $\tilde z = 0$. From (H2), we have that
   \begin{equation}
   \label{expiii1} 1-  m_{KE}(\tilde z, \eta) = a_{2}(\eta) \tilde z^2 + O(\tilde{z}^3),
   \end{equation}
  where by compactness,  $ \inf_{ \mathbb{S}^{d-1}}   |a_{2}(\eta)| \geq c_{s} >0$ for some $c_{s}>0$.
    In particular, we find that
$$ |\xi|^2 \tilde z^2 + 1 -  m_{KE}(\tilde z, \eta) = (a_{2}(\eta) -  |\xi|^2 ) \tilde z^2  + O(\tilde{z}^3)$$
 and hence that for $\eps_{1}$ and hence $|\xi|$ sufficiently small, 
 $$ \left| |\xi|^2 \tilde z^2 + 1 -  m_{KE}(\tilde z, \eta) \right| \geq \frac{c_s}{2} |\tilde{z}|^2.$$
  This also yields that $m_{VP} \over 1-m_{VP}$ is uniformly bounded thanks to 
   \eqref{expiii2}--\eqref{expiii1}.
 
%$$ \left|1-m_{VP}(\gamma,\tau,\xi)\right|   
%=  { 1 \over | \gamma + i \tau |^2}  \left| (\gamma + i \tau)^2  +  (1 - m_{KE}(\gamma, \tau, \xi) \right|.
%$$
%Since $\mu$ satisfies~\eqref{Penrose-kinEuler} and $m_{KE}$  is homogeneous of degree zero we have that
% $$ | 1 - m_{KE}(\gamma, \tau, \xi)| \geq \kappa >0$$
% uniformly for $(\gamma, \tau, \xi) \in C_{R_{1}}$ for some $R_{1}\in (0, R_{0}).$
%  Indeed, by zero homogeneity, it is equivalent to prove the estimate on $\overline{\Omega_{R_{1}} }$ which is compact. 
%   From (H2), the estimate is satisfied if we assume in addition that $\gamma \geq 0$. By continuity it is still true
%   for $\gamma = c | \xi|$ with $ c \in [-R_{1}, 0]$ for some $R_{1}$ sufficiently small.
%   Consequently, for $\eps_{1}$ sufficiently small, we get that
%$$ \left|1-m_{VP}(\gamma,\tau,\xi)\right| \geq  { 1 \over 2 }  { 1 \over | \gamma + i \tau |^2}  \kappa.$$
%This yields iii).  

Next, we prove  iv). We use again that 
\begin{equation}
\label{impl1} 1-m_{VP}(\gamma,\tau,\xi)=   { 1 \over (\gamma + i \tau )^2} \left( ( \gamma + i\tau)^2  +  (1 - m_{KE}(\gamma, \tau, \xi)\right).
\end{equation}
This yields, as $|\tau|^2 \leq 1/\eps_1^2$,
$$  \left | 1-m_{VP}(\gamma,\tau,\xi) \right| \geq  { \eps_{1}^2 \over 2} \left|  ( \gamma + i\tau)^2  +  (1 - m_{KE}(\gamma, \tau, \xi) )\right|.
$$
We shall need the behavior  of  $m_{KE}(\gamma, \tau, \xi)$ close to $\xi=0$. %, $\gamma =0$ and $\tau= \pm 1$.
%   Since the coefficients will be useful later, we give the precise expansion in the next lemma.
%   \begin{lem}
%   \label{lem-symbolkE-sphere}
%We have the expansion  
%\begin{equation}
%m_{KE} (\gamma, \tau, \xi)= { 1 \over (\gamma + i \tau)^2} \left(  \int_{\R^d} \mu(v) |v|^2 \, dv |\xi|^2 +
%|\xi|^2 m_{2}(\gamma, \tau, \xi) \right)
%\end{equation}
%where  $m_{2}(0, \tau,0)= 0$, and for $| \xi | \leq { 1 \over 4}$ and  $(\gamma, \tau, \xi) \in \Omega_{R_{0}/2}$, 
%we have the estimates 
%$$  | \partial_{z}^\alpha \partial_{\xi}^\beta m_{2}(\gamma, \tau, \xi) | \lesssim 
%     1, \qquad  (\gamma, \tau, \xi) \in \Omega_{R_{0}/2}, \, \quad
%      | \xi | \leq { 1 \over 4}.$$
%  \end{lem}
%We postpone the proof of the lemma until the end of the section.
By using the expansion \eqref{mKErefined} for $l=1$, we obtain that  
 in this regime,   for some $C>0$,
$$ |m_{KE}(\gamma, \tau, \xi)| \leq { C \over  \eps_{1}^2} | \xi |^2, $$
therefore, we obtain that for $\eps_{2}$ sufficiently small
 \begin{multline}
 \label{bof} \left|  ( \gamma + i\tau)^2  +  (1 - m_{KE}(\gamma, \tau, \xi) )\right| \geq
  | \gamma + i( \tau -1) | | \gamma + i(\tau + 1) | - | m_{KE}(\gamma, \tau, \xi) |
   \\ \geq  | \tau - 1 | \, |\tau + 1 | - { C \over \eps_{1}^2} | \xi |^2 
 %  \geq   \frac{2}{3} A \eps | \xi |  -   {C \eps  \over \eps_{1}^2} | \xi |
    \geq   \left( \frac{2}{3}  A - {C \over \eps_{1}^2 } \right) \eps | \xi|.
 \end{multline}
We thus find \eqref{belowmiddle} for $\eps_{2}$ sufficiently small and $A\geq A_{0}$ sufficiently large.

We finally prove v). We use again \eqref{impl1}, we have to study the zeroes of
$$ g(z, \xi)=  z^2 + 1 - m_{KE}(\gamma, \tau, \xi).$$
 Writing $z = \pm i +   r \mathfrak{z}$, $\xi= r \omega$, with $| \mathfrak{z}|$ small,  we get by using 
 Lemma \ref{remarkmKE} and the expansion \eqref{mKErefined} that 
 \begin{equation}
\label{defm2}
 g(z,\xi)=  \pm 2 i r  \mathfrak{z} + r^2  \mathfrak{z}^2  -  {1 \over  (\pm i + r  \mathfrak{z})^2} \left( 3 H_{\mu}\omega \cdot \omega\,  r^2  + r^4 m_{2}(\pm i  + r  \mathfrak{z}, r, \omega) \right),
\end{equation}
 where $m_{2}$ is a smooth function of its arguments. We can thus set
 $$ g(z,\xi)
 = r f_{\pm}( \mathfrak{z}, r, \omega),$$
  where 
  \begin{equation}
  \label{f+def}f_{\pm}( \mathfrak{z}, r, \omega) =  \pm 2 i   \mathfrak{z} + r  \mathfrak{z}^2  -  { r \over  (\pm i + r  \mathfrak{z})^2} \left( 3 H_{\mu} \omega \cdot \omega+   r^2 m_{2}(\pm i  + r  \mathfrak{z}, r, \omega)\right).
  \end{equation}
  It thus suffices to study the zeroes of $f_{\pm}$ for $| \mathfrak{z} |$ sufficiently small, $r>0$ close to zero.
   We would like to use the implicit function theorem, nevertheless, since $r=0$ is on the boundary of the domain of definition
    of $f$, we shall first look for a smooth  extension of $f_{\pm}$  for small negative $r$.
     We can use again the expansion \eqref{mKErefined} in Lemma \ref{remarkmKE}  to observe that
    $r^2 m_{2}$ can be expanded as a polynomial in $r$ with even powers plus a remainder of order $\mathcal{O}(r^{2n})$.
     Consequently, we choose an extension by setting 
   $$ m_{\pm}( \mathfrak{z}, r, \omega)=  r^2 {m}_{2} (\pm i  + |r|  \mathfrak{z}, |r| , \omega), $$
then $m_{\pm}$ is a $\mathscr{C}^1$ function %(we identify $\mathbb{C}$ and $\mathbb{R}^2$) 
of its arguments for $| \mathfrak{z}| \leq R_{0}/2$, $|r|<1/2$, $\omega \in \mathbb{S}^{d-1}$, which moreover satisfies
\begin{equation}
\label{flat}
m_{\pm}( \mathfrak{z},0, \omega )  = \partial_{r}m_{\pm}( \mathfrak{z},0, \omega )= 0.
\end{equation}
Let us set
$$ F_{\pm} ( \mathfrak{z}, r, \omega) =  \pm 2 i   \mathfrak{z} + r  \mathfrak{z}^2  -  { r \over  (\pm i + r  \mathfrak{z})^2} \left(  3H_{\mu} \omega \cdot \omega +  m_{\pm}( \mathfrak{z}, r, \omega)\right)$$
and observe that $F_{\pm}$ is a  $\mathscr{C}^1$, $\mathbb{C}$ valued  function of its arguments for $| \mathfrak{z}| \leq R_{0}/2$, $|r|<1/2$, $\omega \in \mathbb{S}^{d-1}$ that coincides with $f_{\pm}$ if $r >0$. Therefore it suffices to study the zeroes of $F_{\pm}$.
For every $\omega \in  \mathbb{S}^{d-1}$,  using  \eqref{flat},    $ F_{\pm} ( 0, 0, \omega)= 0 $ and                                               
  $$D_{w} F _{\pm} ( 0, 0, \omega)= \pm 2 i $$
  is invertible (as a linear map from $\mathbb{R}^2$ to $\mathbb{R}^2$). Therefore by the implicit function theorem, for every $\omega \in   \mathbb{S}^{d-1}$, there exists
    a vicinity of  $(0, 0, \omega)$ such that the zeroes of $F_{\pm}$ are given by a  $\mathscr{C}^1$ curve. 
      By compactness, we can then find $\eps _{3}$
      such that for every  $|r| \leq \eps_{3}$, $ | \mathfrak{z}| \leq \eps_{3},$ and $\omega \in \mathbb{S}^{d-1}$, the zeroes  of $F_\pm(\cdot, r,\omega)$
are described by a curve   $\mathfrak{z} = W_{\pm}(r, \omega)$ such that $W_{\pm}(0, \omega)=0$.
Since by using again \eqref{flat}, we have 
$$\partial_{r}  F_{\pm}( 0, r, \omega)= 3   H_{\mu} \omega \cdot \omega \neq 0, $$
and we also obtain that 
$$ \partial_{r} W_{\pm}(0, \omega) =  \pm i {3 \over 2} H_{\mu} \omega \cdot \omega.$$
This yields v). Note that
 $$ \partial_{r}   \Gamma_{\pm}(0, \omega) = 0, \quad   \partial_{r}   \mathcal{T}_{\pm}(0, \omega) = \pm {  3 \over 2} H_{\mu} \omega \cdot \omega.$$ 
  The fact that we necessarily have  $\Gamma_{\pm}(r, \omega)<0$ for $r>0$ is  a consequence of (H1).
\end{proof}

% This ends the proof of Proposition \ref{lem-zeroes}.  There only remains to prove Lemma \ref{lem-symbolkE-sphere}.
%  \end{proof}
%
%
%
%
%
%\begin{proof}[Proof of Lemma \ref{lem-symbolkE-sphere}]
%
%By integration by parts in $s$
%\begin{align*}
%m_{KE}(\tilde\gamma, \tilde\tau, \tilde\xi)= \frac{1}{\tilde\gamma +i\tilde\tau} m_1 (\tilde\gamma, \tilde\tau, \tilde\xi)
%\end{align*}
%since
%$$
%\int_{\R^d} v_k v_k \na_v \mu \, d v =0,
%$$
%with
%\begin{equation}
%m_1 (\tilde\gamma, \tilde\tau, \tilde\xi) =  - \int_{0}^{+ \infty} e^{-(\gamma + i \tau)s}\, { \tilde\xi \over  |\tilde\xi|^2}\cdot   \sum_{k,l,j} \tilde\xi_k \tilde\xi_l \tilde\xi_j \mathcal{F}_v( v_k v_l v_j \na_v \mu)(\tilde\xi s) \, d s.
%\end{equation}
%This symbol is homogeneous of order $1$.
%Through another integration by parts, we obtain
%\begin{align*}
%m_{KE}(\tilde\gamma, \tilde\tau, \tilde\xi)= \frac{1}{(\tilde\gamma +i\tilde\tau)^2} \left[  \int_{\R^d} \mu(v) |v|^2 \, dv |\tilde\xi|^2    + m_2 (\tilde\gamma, \tilde\tau, \tilde\xi) \right],
%\end{align*}
%since
%$$
%{ \tilde\xi \over  |\tilde\xi|^2}\cdot   \sum_{k,l,j} \tilde\xi_k \tilde\xi_l \tilde\xi_j \mathcal{F}_v( v_k v_l v_j \na_v \mu)(0)  =  -\int_{\R^d} \mu(v) |v|^2 \, dv |\tilde\xi|^2 ,
%$$
%with 
%\begin{equation}
%m_2 (\tilde\gamma, \tilde\tau, \tilde\xi) =   \int_{0}^{+ \infty} e^{-(\gamma + i \tau)s}\, { i \tilde\xi \over  |\tilde\xi|^2}\cdot   \sum_{k,l,j,m} \tilde\xi_k \tilde\xi_l \tilde\xi_j \tilde\xi_m \mathcal{F}_v( v_k v_l v_j  v_m \na_v \mu)(\tilde\xi s) \, d s.
%\end{equation}
%This symbol is of order $2$.

In the above proof, we have used the implicit function theorem in polar coordinates in order to describe the zeroes
 of $z^2 + 1  - m_{KE}$ in the region $|\gamma | \leq \eps_{3} | \xi|$, $ | \tau \pm 1| \leq \eps_3|\xi| $ and $|\xi | \leq \eps_{3}$.
 Nevertheless, it will be useful to get that $Z_{\pm}(r, \omega)$ are actually smooth functions of $\xi$
 under the additional symmetry assumption \eqref{muradial}.
 \begin{lem}
 \label{lemsmooth}
   Assuming  \eqref{eq:norma}, \eqref{muanal}, \eqref{muodd}, \eqref{muradial} and (H1), there exists $\eps_{3}>0$ such that for every $\xi \neq 0$, the zeroes
of $1-m_{VP}$  or equivalently of $ z^2  + 1 - m_{KE}$  with $|\gamma | \leq \eps_{3} | \xi|$, $ | \tau \pm 1 | \leq \eps_3|\xi| $ and $|\xi | \leq \eps_{3}$ are given by two 
  smooth curves of class $\mathscr{C}^{d+2}$  under the form 
     $$ Z_{\pm}(\xi) = \pm i + i |\xi|^2 \Phi_{\pm}(\xi)$$
    where
    $$ \Phi_{\pm}(0) = \pm {3 \over 2} C_{\mu} \in \mathbb{R},  \,   \operatorname{Im} \, \Phi_{\pm} \geq 0.$$
 \end{lem}
 \begin{proof}
 By using the notations of the proof of Proposition \ref{lem-zeroes} v),  since $W_\pm$ is $\mathscr{C}^1$ and $W_{\pm}(0,\omega)= 0$, we can set 
  $W_{\pm}(r, \omega )= r \tilde{W}_{\pm}(r,\omega)$ and thanks to \eqref{f+def}, we see that for $r \neq 0$, 
  $\tilde{W}_{\pm}(r,\omega)$ is a zero of $\tilde f_{\pm}( \mathfrak{z},r,\omega)$
  where
  $$  \tilde f_{\pm}(  \mathfrak{z},r,\omega)   
   =  \pm 2 i    \mathfrak{z} + r^2  \mathfrak{z}^2  -  { 1 \over  (\pm i + r^2  \mathfrak{z})^2} \left( 3 H_{\mu} \omega\cdot \omega +   r^2 m_{2}(\pm i  + r^2  \mathfrak{z}, r, \omega)\right).$$ 
   Moreover, thanks to \eqref{muradial}, we have that 
   $$ H_{\mu} \omega \cdot \omega=   C_{\mu}$$ is independent of $\omega$ and that
    by using Lemma \ref{remarkmKE} and in particular the expansion \eqref{mKErefined} and \eqref{Q2pbetter}, we infer that 
    $ r^4 m_{2}(\pm i  + r^2  \mathfrak{z}, r, \omega) $ has an expansion in terms of polynomials of  $r^2$  of valuation  larger than two plus a high order remainder of the form~\eqref{resteKE}. We can therefore write
    $$ r^2 m_{2}(\pm i  + r^2  \mathfrak{z}, r, \omega) =: m_{\pm}( \mathfrak{z}, \xi)$$
    where $m_{\pm}$ is a smooth function  of its arguments such that $m_{\pm}( \mathfrak{z}, 0) = 0, $ $D_{ \mathfrak{z}}m_{\pm}( \mathfrak{z},0)=0$.
    We can thus 
    write $\tilde f_{\pm}$ as a smooth function of $\xi$:
    $$  \tilde f_{\pm}(  \mathfrak{z},\xi)
    = \pm 2 i    \mathfrak{z} + | \xi|^2  \mathfrak{z}^2  -  { 1 \over  (\pm i + |\xi|^2  \mathfrak{z})^2} \left( 3 C_{\mu} +    m_{\pm}( \mathfrak{z}, \xi)\right).
    $$ 
    Moreover, we observe that
  $$ \tilde f_{\pm}( \pm{3 \over 2} i C_{\mu}, 0)= 0,  
   \quad D_{ \mathfrak{z}}  \tilde f_{\pm}( \pm {3 \over 2}   i C_{\mu}, 0)= \pm 2i.$$
    Consequently, from the implicit function theorem we find that $\tilde{W}_{\pm}(r, \omega)$ is a smooth function of $\xi$
     that we still denote by  $ \tilde{W}_{\pm}(\xi)$. This yields
     $$ Z_{\pm}(r, \omega) = | \xi|^2 \tilde{W}_{\pm}(\xi) , \quad   \tilde{W}_{\pm}(0) = \pm {3 \over 2} i C_{\mu},$$
     which concludes the proof of the lemma.
  \end{proof}  
 
\section{Kernel estimates}

\subsection{Short time estimates}
We start with short time estimates, which require little assumption on $\mu$.
\begin{prop}
\label{lemshort}
Assuming~\eqref{muanal}, there exists $C>0$ such that for every $t \in(0, 1]$, 
$$ \|G(t) \|_{L^1} \leq C t, \quad   \| G(t) \|_{L^\infty} \leq C { 1 \over t^{d-1}}.$$
\end{prop}

\begin{proof}
We observe that $ G(t,x)$ solves the integral equation
$$   G=  K  + K *_{t,x}  G$$
\begin{align*} K(t,x)= e^{\gamma t} \mathcal{F}^{-1}( m_{VP}(\gamma, \tau, \xi) ) (t,x) & = \int_{\mathbb{R}^d \times \mathbb{R}}
 e^{ \gamma t  + i \tau t  + i x \cdot \xi} m_{VP}(\gamma, \tau, \xi)\, d\tau d\xi \\
 &  = - { 1 \over t^{d-1}} \mu\left(  { x \over t} \right)
  \mathrm{1}_{t \geq 0}.
  \end{align*}
  Therefore,  we have the estimate
  $$ \|K(t)\|_{L^1} \lesssim t.$$
  Moreover, as already observed, $G$ vanishes in the past, therefore
  $$ G(t,x) =  K(t,x) + \int_{0}^t K(t-s, \cdot)*_{x} G(s, \cdot) \, ds, \quad  \forall t \geq 0.$$
   This yields
   $$ \|G(t) \|_{L^1} \lesssim  t +  \int_{0}^t (t-s)  \|G(s)\|_{L^1}\, ds$$ and hence from the Gronwall
   inequality, we get that
   $$ \|G(t) \|_{L^1} \lesssim t, \quad \forall t \leq 1.$$  
  We also obtain  that
   $$ \|G(t) \|_{L^\infty} \lesssim { 1 \over t^{d-1}} + \int_{0}^{{t\over 2} } { 1 \over (t-s)^{d-1}} \|G(s) \|_{L^1} \; ds+ \int_{t\over 2}^t
     (t-s) \|G(s) \|_{L^\infty} \; ds.$$
    This yields for $t \leq 1$ that $y(t)= t^{d-1} \|G(t)\|_{L^\infty}$ verifies
    $$ y(t) \lesssim 1 +  t^{d} \int_{0}^{t\over 2} { 1 \over (t-s)^{d-1}} \, ds + t^{d-1}  \sup_{[0, t]} y(s) 
     \int_{t\over2}^t  (t-s)  { 1 \over s^{d-1}} \, ds.$$
     Since 
     $$  t^{d} \int_{0}^{t\over 2} { 1 \over (t-s)^{d-1}} ds = t^2 \int_{0}^{ 1 \over 2} { 1 \over (1- u)^{d-1}} \, ds = Ct^2$$
      and 
      $$    t^{d-1}    \int_{t\over2}^t  (t-s)  { 1 \over s^{d-1}} \, ds \leq  2^{d-1} t^2 \int_{{1 \over 2}}^1 (1-u)\, du \lesssim  t^2, $$
      we get that for every $T >0$
      $$ \sup_{[0, T]} y(t) \lesssim 1 + T^2 + T^2   \sup_{[0, T]} y(t).$$
       This yields the result for $t\in (0, T]$, $T$  sufficiently small. We can then iterate the argument finitely many times in a classical way  to get the result for $t \in (0, 1].$
       This ends the proof.
   \end{proof}
   
   \subsection{Large time estimates}

   We shall now focus on estimates for $t \geq 1$. First, observe that by setting $z= \gamma + i \tau$,  we can write \eqref{Gdef} as
$$ G (t,x)={ 1 \over i} \int_{\mathbb{R}^d} e^{i x \cdot \xi} \left( \int_{\operatorname{Re}\, z= \gamma }   e^{z t } { m_{VP}(z,\xi) \over 1-m_{VP}(z, \xi)}
 \, dz \right) d \xi.$$
 Let us pick $\delta>0$ to be fixed later.
   We split $G$ as  a high frequency and a low frequency part:
   \begin{equation}
   \label{Gsplit1}
    G (t,x)= G^H(t,x)+ G^L(t,x)
   \end{equation}
   where
   \begin{align}
  \label{GHdef}  G^H(t,x) & = { 1 \over i} \int_{\mathbb{R}^d} e^{i x \cdot \xi} \left( \int_{\operatorname{Re}\, z= \gamma }   e^{z t } { m_{VP}(z,\xi) \over 1-m_{VP}(z, \xi)}
   \left(1- \chi\left({\xi \over \delta}\right) \right)
 \, dz \right) d \xi, \quad \\
 \label{GLdef} G^L(t,x) & = { 1 \over i} \int_{\mathbb{R}^d} e^{i x \cdot \xi} \left( \int_{\operatorname{Re}\, z= \gamma }   e^{z t } { m_{VP}(z,\xi) \over 1-m_{VP}(z, \xi)}
   \chi\left({\xi \over \delta} \right) \, dz \right) d \xi.
   \end{align}
   where $\chi \in \mathcal{C}^\infty_{c}(\mathbb{R}^d)$ is a nonnegative radial  function equal to one for $| \xi| \leq 1$
    and supported in the  ball of radius $2$. Note that $G^H$ and $G^L$ depends on $\delta$. The choice of $\delta$
     will be carefully performed in order to estimate $G^L$.

\subsubsection{High frequency estimates}
 
    We shall first estimate the high frequency contribution $G^H$. 
    %, {\color{red}which like the torus case \cite{GNR} gives an exponential decay in time. }

  \begin{prop}
  \label{lemGH}
  Assuming \eqref{eq:norma}, \eqref{muanal}, \eqref{muodd} and (H1), 
  for every $\delta>0$, 
   there exist $C>0$ and $\alpha>0$ such that 
   $$ \|G^H(t) \|_{L^1} \leq C e^{- \alpha t}, \quad  \|G^H(t) \|_{L^\infty} \leq  C e^{- \alpha t}, \quad \forall t \geq 1.$$
  \end{prop}
 
 \begin{proof}
 Let us first recall that  $(1 - \chi)$ is supported in the zone $|\xi| \geq \delta >0$ so that the argument is very similar
 to the one used in the torus case in \cite{GNR} for example.
  For $\xi \neq 0$, thanks to the Penrose stability condition (H1) and Corollary \ref{lem-symbolVP-sphere}, the function    $ m_{VP}(\cdot, \xi)/(1-m_{VP}(\cdot, \xi))$ is an holomorphic
   function in $\{\operatorname{Re} z >0\}$. Moreover, by using  ii) of Proposition \ref{lem-zeroes}, for $|\xi| \geq \delta$,  it extends as an holomorphic function
    in  $ \{\operatorname{Re} z \, >  - R_{\delta} | \xi|\} $,   where $R_{\delta}$ is given by ii) of Proposition \ref{lem-zeroes}, and we have a positive uniform estimates from below of  $|1-m_{VP}|$. We can then 
    use the Cauchy formula to get that
    for $|\xi| \geq \delta$, 
$$ \int_{\operatorname{Re}\, z= \gamma }   e^{z t } { m_{VP}(z,\xi) \over 1-m_{VP}(z, \xi)}
   \chi\left({\xi \over \delta} \right) \, dz  = \int_{\operatorname{Re}\, z= - R_{\delta} |\xi|/2}   e^{z t } { m_{VP}(z,\xi) \over 1-m_{VP}(z, \xi)}
   \chi\left({\xi \over \delta} \right) \, dz. 
   $$
    Indeed, we can apply  Corollary \eqref{lem-symbolVP-sphere} to get that uniformly  for  $|\xi| \geq \delta$ and $\operatorname{Re} z \geq -R_{\delta} |\xi|/2$, 
    \begin{equation}
    \label{mVPhigh}
     |m_{VP}(z, \xi) | \lesssim {  C \over |\xi|^2  + \tau^2}
    \end{equation}
    so that there is no contribution from infinity.
    Consequently, we have to estimate
    $$   G^H(t,x)=   \int_{\mathbb{R}^d } e^{i x \cdot \xi} \int_{\mathbb{R}}   e^{- R_{\delta} | \xi| t/2 }  e^{i \tau t}{ m_{VP}(- R_{\delta} | \xi| /2,\tau,\xi) \over 1-m_{VP}(- R_{\delta} | \xi| /2, \tau, \xi)}
   \left( 1- \chi\left({\xi \over \delta}\right)\right)  d\tau d\xi.$$
  By using again \eqref{mVPhigh} and ii) of Proposition \ref{lem-zeroes}, we easily get that
 $$  |G^H(t,x)| \lesssim  \int_{|\xi| \geq \delta }   e^{- R_{\delta} | \xi| t/2 }  \int_{ \mathbb{R}}  {  1 \over |\xi|^2  + \tau^2}
  \, d\tau\, d\xi  \lesssim    \int_{|\xi| \geq \delta }   e^{- R_{\delta} | \xi| t/2 }  { 1 \over | \xi |} \, d \xi
   \lesssim  e^{- \tilde \alpha t} $$
   for some $ \tilde \alpha >0$.
   This yields 
   $$ \| G^H(t) \|_{L^\infty} \lesssim e^{- \tilde \alpha t}, \, \forall t \geq 1.$$
   For the $L^1$ norm,  by integrating by parts in $\xi$ and applying Corollary~\ref{lem-symbolVP-sphere}, we obtain in a similar way that for all  multi-indices $|\beta|\leq d+1$,
   $$ | x^\beta G^H(t,x) | \lesssim (1+   t^{|\beta|}) e^{ - \tilde  \alpha t}, \quad \forall t \geq 1.$$
   Therefore, we obtain that
   $$ | G^H(t,x) | \leq { 1 \over 1 + |x|^{d+1}} (1 +   t^{d+1}) e^{ - \tilde  \alpha t}, \quad \forall t \geq 1$$
    and hence that
    $$ \|G^H(t)\|_{L^1} \lesssim e^{- \alpha t}, \quad \forall t \geq 1,$$
    with $\alpha = \tilde \alpha/2.$
     \end{proof}
  
  \subsubsection{Low frequency estimates}
     
     We shall now estimate the low frequency part $G^L$.
     
  \begin{lem}
  \label{lemGL}
  Assuming \eqref{eq:norma}, \eqref{muanal}, \eqref{muodd}, (H1), (H2), for $\delta >0$ small enough,  
  we have the following  decomposition of $G^L$:
 \begin{equation}
 \label{GLdec-prop}
 G^L(t,x) = G^r(t,x)  + G^S_{+}(t,x) + G^S_{-}(t,x)
 \end{equation}
 where, 
 \begin{align*}
  & G^r(t,x) =   \int_{\mathbb{R}^d} e^{ix \cdot \xi} \left( \int_{\operatorname{Re}\, z=- \tilde \delta |\xi| }   e^{z t } { m_{VP}(z,\xi) \over 1-m_{VP}(z, \xi)}dz \right) \chi \left( { \xi \over \delta}  \right)  \, d\xi, \quad \tilde\delta = \delta^{3/2},\\
 & G^S_{\pm}(t,x) =  2\pi  \int_{\mathbb{R}^d} e^{ix \cdot \xi} e^{ Z_{\pm}(r, \omega) t}   a_{\pm}(\xi)  \chi \left( { \xi \over \delta} \right)  \, d\xi,
  \quad a_{\pm}(\xi)=  { Z_{\pm}(r, \omega)^2 \over  2 Z_{\pm}(r, \omega) - \partial_{z}m_{KE}(Z_{\pm}(r, \omega), \xi)}.
 \end{align*}

  \end{lem}
    \begin{proof}
     We now deal with the region $| \xi| \leq \delta$ for $\delta>0$ to be chosen sufficiently small.     
      For $\xi\neq 0$, we would like again to use the Cauchy formula to change the integration contour for 
     $$  I_{\xi}=  { 1 \over i } \int_{\operatorname{Re}\, z= \gamma }   e^{z t } { m_{VP}(z,\xi) \over 1-m_{VP}(z, \xi)}
   \chi \left({\xi \over \delta} \right) \, dz .$$
   Again for $\gamma >0$ the function  $1-m_{VP}(z, \xi)$ does not vanish thanks to (H1) so that
   we have to carefully study what happens for negative $\gamma$ with  $|\gamma|$ small. 
    We observe that  thanks to i) and iii) of Proposition \ref{lem-zeroes}, for $|(\gamma, \tau, \xi)| \leq 2 \eps_{1}$ or  $|(\gamma, \tau, \xi)| \geq 1/ (2\eps_{1})$ (reducing
   $\eps_{1}$ if necessary), the function
    $1-m_{VP}$ does not vanish in $\mathcal{C}_{R_{1}}$.
    
    We shall now choose
   \begin{equation}
   \label{choixdelta} 
  0 < \delta  \leq  {\eps_{3} \over  10 A_0}
   \end{equation}
   where  $\eps_{3}$ and $A_0$ are given by Proposition \ref{lem-zeroes} iv) and v). As a consequence, for $| \xi | \leq \delta$, $|\gamma| \leq \delta  | \xi| $,  if $| \tau| \geq 1/\eps_{1}$ or $|\tau| \leq \eps_{1}$, 
    $1-m_{VP}$ does not vanish.
   Then, since $\delta \leq \eps_{3}$, we get that 
    for  $| \xi | \leq \delta $,   $|\gamma| \leq \delta |\xi|$ and $|\tau \pm 1 | \leq  \eps_3 | \xi|, $
    the function  $  1-m_{VP} $ has for  each $\xi$  exactly two zeroes  described by  v) of Proposition
   \ref{lem-zeroes}. Moreover, since $\Gamma_{\pm}(0, \omega)= 0$, $\partial_{r} \Gamma_{\pm}(0, \omega)= 0$, 
     we  have that for  $|\tau \pm 1 | \leq \delta  | \xi|$,  the zeroes in  $-\delta | \xi |  \leq \gamma \leq 0$
     are actually localized in   $-\ C\delta | \xi |^2 \leq - C |\xi|^3 \leq \gamma \leq 0$ for some $C>0$.
      Therefore, assuming that $\delta$ is sufficiently small, we get in particular that on the line 
        $\operatorname{Re} z =  \gamma = - \delta^{3\over 2}|\xi|$,  there is no zero of $ 1 - m_{VP}$ for   $|\tau \pm 1| \leq \delta | \xi|$, $|\xi| \leq \delta$.
        Next, using   iv) of Proposition \ref{lem-zeroes}  since  $\delta \leq \eps_{2}$, we get that
        for  $ | \gamma | \leq \delta^{3 \over 2} | \xi|, $  $|\tau \pm 1 | \geq  \delta | \xi| $, $\eps_{1} \leq |\tau| \leq 1/\eps_{1}$ and $| \xi | \leq \delta$, 
        $1-m_{VP}$ does not vanish.

        To summarize, we have thus obtained that for each $\xi \neq 0$,   $|\xi| \leq \delta $, 
        there are exactly two zeroes of $(1-m_{VP})$ in the region
         $|\gamma | \leq \delta |\xi| $ and they are  described by v) of Proposition \ref{lem-zeroes}. Moreover, 
          they are localized in  $-\ C\delta | \xi |^2  \leq \gamma \leq 0$ and $ |\tau \pm 1 | \leq \eps_{3} | \xi|.$
          We can thus use the residue formula to write that for $\xi \neq 0$ (note that there is again no contribution from infinity
           since the estimate \eqref{mVPhigh} is still valid for large $\tau$),
     $$  I_{\xi}= {1 \over i } \int_{\operatorname{Re}\, z=- \delta^{3 \over 2} |\xi| }   e^{z t } { m_{VP}(z,\xi) \over 1-m_{VP}(z, \xi)}
   \chi\left({\xi \over \delta} \right) \, dz  +  2 \pi  \chi\left({\xi \over \delta} \right)  \sum_{\pm} e^{ Z_{\pm}(r, \omega) t}\left(\mbox{Res} { m_{VP} \over 1 - m_{VP}}(\cdot, \xi) \right)_{|Z_{\pm}(r, \omega)}$$
   where $r = |\xi|,$ $\omega = \xi/|\xi|$.       
 Computing the residue, we obtain
 $$ 2 \pi  \sum_{\pm}\left(\mbox{Res} { m_{VP} \over 1 - m_{VP}}(\cdot, \xi) \right)_{|Z_{\pm}(r, \omega)}
 =  - 2 \pi   \sum_{\pm}  { 1 \over \partial_{z} m_{VP}( Z_{\pm}(\xi), \xi)}.$$
    To get regularity in $\xi$ close to $\xi= 0$, it is convenient to express the residue in terms of $m_{KE}$. Thanks to \eqref{eq-K1def}, 
    we have
    $$ m_{VP}= - { 1 \over z^2} ( 1 - m_{KE})$$
    and hence
    $$ \partial_{z} m_{VP} = { 2 \over z^3}   ( 1 - m_{KE})  + { 1 \over z^2} \partial_{z} m_{KE}.$$
     Since $ ( 1 -m_{KE}) = - z^2$ at $z= Z_{\pm}$, we can also write
    $$ 2  \pi  \sum_{\pm}\left(\mbox{Res} { m_{VP} \over 1 - m_{VP}}(\cdot, \xi) \right)_{|Z_{\pm}(r, \omega)}
 =   2 \pi   \sum_{\pm}  { Z_{\pm}(r, \omega)^2 \over  2 Z_{\pm}(r, \omega) - \partial_{z}m_{KE}(Z_{\pm}(r, \omega), \xi)}.$$
 This yields the decomposition of $G^L$
 \begin{equation}
 \label{GLdec}
 G^L(t,x) = G^r(t,x)  + G^S_{+}(t,x) + G^S_{-}(t,x)
 \end{equation}
 where (setting $\tilde \delta = \delta^{3 \over 2}$ for notational convenience), 
 \begin{align*}
  & G^r(t,x) =   \int_{\mathbb{R}^d} e^{ix \cdot \xi} \left( \int_{\operatorname{Re}\, z=- \tilde \delta |\xi| }   e^{z t } { m_{VP}(z,\xi) \over 1-m_{VP}(z, \xi)}dz \right) \chi \left( { \xi \over \delta}  \right)  \, d\xi, \\
 & G^S_{\pm}(t,x) = 2\pi   \int_{\mathbb{R}^d} e^{ix \cdot \xi} e^{ Z_{\pm}(r, \omega) t}   a_{\pm}(\xi)  \chi \left( { \xi \over \delta} \right)  \, d\xi,
  \quad a_{\pm}(\xi)=  { Z_{\pm}(r, \omega)^2 \over  2 Z_{\pm}(r, \omega) - \partial_{z}m_{KE}(Z_{\pm}(r, \omega), \xi)},
 \end{align*}
giving the lemma.   \end{proof}

  \subsubsection{Low frequency estimates: regular part}

 The next step is to estimate $G^r$ and $G^S_{\pm}$ in \eqref{GLdec}. We start with $G^r$. 

 \begin{prop}[ Study of $G^r$]
 \label{propGr}
 Assuming \eqref{eq:norma}, \eqref{muanal}, \eqref{muodd} and (H2), $\delta$ can be chosen small enough so that 
     uniformly for $t \geq 1$
    $$ \| G^r \|_{L^\infty} \lesssim {1 \over t^{d+1}}, \quad   \| G^r \|_{L^1} \lesssim {1 \over t }.$$
 
 \end{prop}

 \begin{proof}
 Let us write that 
$$ G^r(t,x) =  \frac{1}{i} \int_{\mathbb{R}^d} e^{ix \cdot \xi} \left( \int_{\operatorname{Re}\, z = - \tilde \delta  | \xi| }  e^{- \tilde \delta | \xi | t }  e^{i \tau t }{ m_{VP}(z ,\xi) \over 1-m_{VP}(z, \xi)}\, dz \right) \chi \left( { \xi \over \delta} \right)  \, d\xi
 =   \int_{\mathbb{R}^d} e^{ix \cdot \xi}  \tilde{I}_{\xi}  \chi \left( { \xi \over \delta} \right)  \, d\xi.
$$
Thanks to  Proposition \ref{lem-VP-VB},  we observe that
 \begin{multline*} { m_{VP} (z, \xi)  \over 1 - m_{VP}(z,\xi)} = - { 1 \over z^2} { 1 - m_{KE}(z, \xi) \over  1 + { 1 \over z^2}( 1- m_{KE}(z,\xi))}
  =  - {1 -  m_{KE}(z, \xi)  \over z^2  + 1 - m_{KE}(z, \xi) } \\
   = - 1 +  { z^2 \over z^2+ 1 - m_{KE}} = - 1 + { z^2 \over z^2 + 1 } { 1 \over (1 - {m_{KE}\over 1+ z^2})}
    = - 1 + { z^2 \over z^2 + 1 }  + { z^2 \over( z^2 + 1)^2 } { m_{KE} \over (1 - {m_{KE}\over 1+ z^2}) } \\
     = - { 1 \over z^2 + 1} + {z^2 \over z^2+1}   { m_{KE} \over (z^2+ 1 - m_{KE}) }.
  \end{multline*}
   We can thus write
   $$ \tilde{I}_{\xi}= -{1 \over i } \int_{\operatorname{Re}\, z = - \tilde \delta | \xi| } e^{ zt } { 1 \over  1+ z^2}\, dz +  J_{\xi}, \quad J_{\xi}=
   { 1 \over i } \int_{\operatorname{Re}\, z = - \tilde \delta  | \xi| }e^{zt}   {z^2 \over z^2+1}   { m_{KE} \over (z^2+ 1 - m_{KE})}(z, \xi)\, dz.$$
    Next, we observe that from the Cauchy formula
    $$  \int_{\operatorname{Re}\, z = - \tilde \delta | \xi| } e^{ zt } { 1 \over  1+ z^2} \, dz
     =  \int_{\operatorname{Re}\, z = - \Gamma } e^{ zt } { 1 \over  1+ z^2}dz $$
     for any $\Gamma \geq \tilde \delta | \xi|$
     and thus, sending $\Gamma$ to $+ \infty$, we get that
     $$  \int_{\operatorname{Re}\, z = - \tilde \delta | \xi| } e^{ zt } { 1 \over  1+ z^2} \, dz = 0.$$
     We have thus obtained that
   \begin{equation}
   \label{Grbis}
    G^r(t,x)=  \int_{\mathbb{R}^d}  e^{ i x \cdot \xi}J_{\xi} \chi \left( { \xi \over \delta} \right) \, d \xi, \quad
    J_{\xi}= {1 \over i }
    \int_{\operatorname{Re}\, z = - \tilde \delta  | \xi| }e^{zt}   {z^2 \over z^2+1}   { m_{KE} \over (z^2+ 1 - m_{KE})}(z, \xi)\, dz. 
   \end{equation} 
 We shall use this more convenient form to prove the estimates.
We shall further  split the $J_{\xi}$ term  into
\begin{multline}
\label{splitIxi1}
J_{\xi}=  \int_{| \tau| \leq \eps_{1} }   e^{- \tilde \delta | \xi | t }  e^{i \tau t }{z^2 \over z^2+1}   { m_{KE} \over (z^2+ 1 - m_{KE})}(z, \xi)\, d\tau
 \\+  \int_{| \tau| \geq  \eps_{1}, \, ||\tau|- 1 | \geq 1/2  }   e^{- \tilde \delta | \xi | t }  e^{i \tau t }  {z^2 \over z^2+1}   { m_{KE} \over (z^2+ 1 - m_{KE})}(z, \xi)  \, d\tau
  \\ +    \int_{| \tau| \geq  \eps_{1}, \, ||\tau|- 1 | \leq  { 1 \over 2} }   e^{- \tilde \delta | \xi | t }  e^{i \tau t }
 {z^2 \over z^2+1} { m_{KE} \over (z^2+ 1 - m_{KE})}(z, \xi)
  \, d\tau
 =: J_{\xi, 1} + J_{\xi, 2} + J_{\xi, 3},
\end{multline}
where  we recall that $\eps_{1}$ is defined in iii) of Proposition \ref{lem-zeroes} 
and we decompose accordingly  $G^{r} $ into
\begin{equation}
\label{splitGr} G^r =  G^r_{1}+ G^r_{2} + G^r_{3}.
\end{equation}
In the following three lemmas, we shall provide estimates of $G^r_{i}$, $i=1, \, 2, \, 3$.

Let us start with $G^r_{2}$, which is the easiest one.
\begin{lem}
\label{lemGr2} Under the assumptions of Proposition \ref{propGr}, 
 $\delta$ can be chosen small enough so that  uniformly for $t \geq 1$,
 $$ \|G^r_{2}\|_{L^\infty} \lesssim { 1 \over t^{d+1}}, \quad   \|G^r_{2}\|_{L^1} \lesssim { 1 \over t}.$$
\end{lem}
\begin{proof}
 By using the same factorization as in  \eqref{bof}, we observe that  
uniformly   for $ |\xi| \leq \delta$, $| |\tau|-1| \geq 1/2$ and $\gamma= -\tilde \delta |\xi|$,  we can take $\delta$ small enough so that
 \begin{equation}
 \label{below1} | z^2 + 1 -m_{KE}( - \tilde \delta | \xi|, \tau, \xi) | \geq \kappa_{0}  |z^2 + 1|>0,
 \end{equation} 
where $\kappa_{0}$ depends only on $\delta$ and $\eps_{1}$.  
Moreover, still in the same range of parameters, 
\begin{equation}
\label{below2} |z^2+ 1| = |z+i| |z- i| \gtrsim  1 + \tau^2, \quad |z|^2 \lesssim |\xi|^2 + \tau^2.
\end{equation}
  Therefore, we get that 
  \begin{multline*}
  |G^r_{2}(t,x) | \lesssim   \int_{| \xi | \leq \delta }  e^{ - \tilde \delta | \xi | t} \int_{\mathbb{R}} {   |\xi|^2 + |\tau |^2 \over 1 + |\tau |^4} |m_{KE}(-  \tilde \delta | \xi |, \tau, \xi) | \, d\tau d \xi
 \\ \lesssim 
   \int_{| \xi | \leq \delta }  e^{ - \tilde \delta | \xi | t} \int_{\mathbb{R}} |m_{KE}(-  \tilde \delta | \xi |, \tau, \xi) | \, d\tau d \xi.
   \end{multline*}
   We then set $\tau= | \xi | \tau'$ and use the degree zero  homogeneity of $m_{KE}$ to get
   $$  |G^r_{2}(t,x) | \lesssim   \int_{| \xi | \leq \delta}  e^{ - \tilde \delta | \xi | t} | \xi|  \int_{\mathbb{R}}
   \left|m_{KE}\left(-  \tilde \delta , \tau', {\xi \over | \xi|}\right) \right| \, d\tau' d \xi. $$
   From  Proposition \ref{lem-VP-VB}, we have that  $\left|m_{KE}\left(-  \tilde \delta , \tau', {\xi \over | \xi|}\right)\right|$
    is uniformly bounded for $|\tau'| \leq 2$, while we get from Lemma \ref{remarkmKE}
     $$   \left|m_{KE}\left(-  \tilde \delta , \tau', {\xi \over | \xi|}\right) \right| \lesssim { 1 \over (\tau')^2}, \quad
      |\tau'| \geq 2,$$
       and hence 
       \begin{equation}
       \label{estimmke}
   \left|m_{KE}\left(-  \tilde \delta , \tau', {\xi \over | \xi|}\right) \right| \lesssim { 1 \over 1 +  (\tau')^2}.
\end{equation}
This yields
   $$  |G^r_{2}(t,x) | \lesssim   \int_{| \xi | \leq \delta}  e^{ - \tilde \delta | \xi | t} | \xi|  \int_{\mathbb{R}}  { 1 \over 1 +  (\tau')^2}\, d\tau' d \xi. $$
   and hence by finally setting $\tilde \xi = t \xi$, we obtain
   $$ |G^r_{2}(t,x) | \lesssim { 1 \over t^{d+1}}.$$ 
   
   There remains to estimate the $L^1$ norm.
    To this end, we use an homogeneous Littlewood-Paley decomposition. We write 
   \begin{multline}
   \label{G2rLP} G_{2}^r = \sum_{q \leq 0} G^r_{2,q}, \quad
    G^r_{2, q}(t,x)= \int_{\mathbb{R}^d}  e^{ i x \cdot \xi}J_{\xi, 2} (t) \chi \left( { \xi \over \delta} \right)   \phi \left( { \xi \over 2^q} \right)\, d \xi,  \\
    J_{\xi, 2}(t)= 
    \int_{|\tau | \geq \eps_{1}, \,  || \tau | - 1| \geq 1/2 }e^{-  \tilde \delta |\xi| t}  e^{i \tau t}  {z^2 \over z^2+1}   { m_{KE} \over (z^2+ 1 - m_{KE})}(- \tilde \delta | \xi|, \tau, \xi)\, d\tau, \quad z= - \tilde \delta |\xi|  + i \tau,
    \end{multline}
    where $\phi$ is supported in the annulus $1/4 \leq | \xi | \leq 4$.
    Changing $\xi$  for  $\xi/ 2^q$, we get that
\begin{multline}
\label{G2-G2-0}
 G^r_{2, q}(t,x) = 2^{qd} \mathcal{G}^r_{2,q} (T, X), \quad \mathcal{G}^r_{2,q} (T, X)
  = \int_{\mathbb{R}^d}  e^{i X \cdot \xi} J_{\xi, 2,q} (T)   \chi \left( { 2^q  \xi \over \delta} \right)   \phi \left(  \xi \right)\, d \xi, 
  \quad  T= 2^qt, \, X= 2^q x,\\
   J_{\xi, 2, q}(T)=
    \int_{ |\tau | \geq \eps_{1}, \, || \tau | - 1| \geq 1/2 }e^{-  \tilde \delta |\xi| T}  e^{i \tau T \over 2^q}  {z^2 \over z^2+1}   { m_{KE} \over (z^2+ 1 - m_{KE})}(
    - \tilde \delta 2^q | \xi |,  \tau,  
    2^q  \xi )\, d\tau ,
    \end{multline} 
    where now the integral in $\xi$ is supported on the annulus $ 1/4 \leq | \xi | \leq 4.$
    We then observe that since $2^q \leq 1$, we have that 
    $$ \| \mathcal{G}_{q} (T, \cdot)\|_{L^1}
     \lesssim  \sum_{ | \alpha | \leq d+1}   \| \partial_{\xi}^\alpha J_{\xi, 2, q}\|_{L^\infty}.$$
      By using again \eqref{below1}, \eqref{below2} and Proposition \ref{lem-VP-VB}, we get
       that  for $|\alpha|\leq d+1$, and uniformly for $1/4 \leq | \xi | \leq 4$, 
       $$    | \partial_{\xi}^\alpha J_{\xi, 2, q}(T) | \lesssim  \int_{ || \tau | - 1| \geq 1/2 }e^{-  \tilde \delta |\xi| T} 
         \sum_{k \leq | \alpha |} 2^{qk} | h_{k} ( - \tilde \delta 2^q | \xi |,  \tau,  
    2^q  \xi )  |\, d\tau,
    $$
    where $h_{k}$ is  positively homogeneous of degree $-k$. Moreover, by using the expansion  of $m_{KE}$ provided by 
    Lemma \ref{remarkmKE}, we also know that 
     for bounded $ \xi$ and $|\tau| \geq 1/2$, 
   $$| h_{k}(z, \xi) | \lesssim  { 1 \over |z^2 +  1 |} \lesssim { 1 \over 1 + \tau^2} .$$
         Therefore, by setting $\tau= 2^q | \xi|  \tau'$, we obtain
     that
   \begin{multline*}  | \partial_{\xi}^\alpha J_{\xi, 2, q}(T) | \lesssim  2^q  \int_{ | 2^q | \xi| \tau' | - 1| \geq 1/2 }e^{- \tilde \delta |\xi| T} 
         \sum_{k \leq | \alpha |} | h_{k} ( - \tilde \delta,  \tau',  
     {\xi \over |\xi| } )  |\, d\tau' \\  \lesssim   2^q  e^{-  \tilde \delta T/2}   \int_{\mathbb{R}} { 1 \over 1 + (\tau')^2 } \, d\tau'
      \lesssim  2^q  e^{-  \tilde \delta T/2}.
      \end{multline*}
 Consequently we obtain from \eqref{G2-G2-0} that
 $$ \|G^r_{2,q}  (t) \|_{L^1} \lesssim  2^q  e^{-  \tilde \delta 2^q t/2} \lesssim { 2^q \over 1 + (2^qt)^4}.$$
  This finally yields for $t \geq 1$
$$
\|G^r_{2} (t) \|_{L^1} \leq \sum_{q \leq 0} \| G^r_{2,q}  (t) \|_{L^1}  \lesssim \sum_{ 2^q \leq 1/t}  2^q
 + {1 \over t^4} \sum_{ 1/t \leq 2^q \leq 0}  { 1 \over 2^{3q}}
  \lesssim { 1 \over t}+ {t^3 \over t^4} \lesssim {1 \over t}.
  $$
   This ends the proof. 
   \end{proof}  
    Let us turn to $G^r_{1}$.
   \begin{lem}
   Under the assumptions of Proposition \ref{propGr},     we have uniformly for $t \geq 1$
    $$ \| G^r_{1} \|_{L^\infty} \lesssim {1 \over t^{d+2}}, \quad   \| G^r_{1} \|_{L^1} \lesssim {1 \over t^{2}}.$$
   \label{lemGr1}
   \end{lem}

   \begin{proof}
    In this regime of low frequencies for $ |\xi| \leq \delta, | \tau | \leq \eps_{1}$, and $\gamma= -\tilde \delta |\xi|$, we have  that
     $$\left| { z^2 \over z^2 + 1} \, { m_{KE} \over z^2 + 1 - m_{KE}}(- \tilde \delta | \xi|, \tau,  \xi) \right|
      \lesssim |z|^2  \left | {m_{KE}  (- \tilde \delta | \xi|, \tau , \xi) \over z^2 + 1 - m_{KE}(  -\tilde \delta | \xi|, \tau,  \xi)} \right| .$$
   By setting again $\tau = |\xi| \tau'$, and by using  that $m_{KE}$ is homogeneous of degree zero,  we obtain that
   $$ |G^r_{1}(t,x)| \lesssim   \int_{| \xi | \leq \delta}  e^{ - \tilde \delta | \xi | t} |\xi|^3   \int_{|\tau' | \leq \eps_{1}/|\xi|} 
    |z'|^2  \left| {m_{KE}  ( - \tilde \delta, \tau' , {\xi \over |\xi|}) \over  |\xi|^2 (z')^2  +  1 - m_{KE}(  -\tilde \delta , \tau ',  {\xi \over |\xi|})}  \right|\, d\tau'$$
    where $z' = -\tilde \delta + i \tau'.$
%    This yields
%    $$  |G^r_{1}(t,x)| \lesssim   \int_{| \xi | \leq \delta}  e^{ - \tilde \delta | \xi | t}   \int_{|\tau | \leq \eps_{1}} 
%    (| \xi|^2 + \tau^2)  \left |m_{KE}  ( \tilde \delta | \xi|, \tau , \xi) \right|\, d\tau.$$
%     By setting again $\tau = |\xi| \tau'$, and that $m_{KE}$ is homogeneous degree zero,  we obtain that
%     $$   |G^r_{1}(t,x)| \lesssim   \int_{| \xi | \leq \delta}  e^{ - \tilde \delta | \xi | t} |\xi|^3   \int_{|\tau' | \leq \eps_{1}/|\xi|} 
%    ( 1  + (\tau')^2)  \left |m_{KE}  ( \tilde \delta, \tau' , {\xi \over |\xi|}) \right|\, d\tau'.$$
By using Proposition \ref{lem-zeroes} iii), in particular \eqref{penrosepetit}, we know that
$$ \left| |\xi|^2 (z')^2  +  1 - m_{KE}(  -\tilde \delta , \tau ',  {\xi \over |\xi|}) \right|$$
 is bounded from below by a positive constant  since $|z'| \geq  \tilde \delta$
and hence, obtain that
$$ |G^r_{1}(t,x)| \lesssim   \int_{| \xi | \leq \delta}  e^{ - \tilde \delta | \xi | t} |\xi|^3   \int_{|\tau' | \leq \eps_{1}/|\xi|} 
   |z'|^2  \left| m_{KE}  ( - \tilde \delta, \tau' , {\xi \over |\xi|})   \right|\, d\tau'.$$
   % Next, we have from Proposition \ref{lem-VP-VB} that $ \left|m_{KE}  \left( \tilde \delta, \tau' , {\xi \over |\xi|} \right) \right|$ is bounded for $|\tau'| \leq 2$, and from 
  %  the expansion of Lemma \ref{remarkmKE}
  %  that
  As in~\eqref{estimmke}, we have
   $$  \left|m_{KE}  \left( - \tilde \delta, \tau' , {\xi \over |\xi|} \right) \right|\lesssim { 1 \over 1+ (\tau')^2}.$$
   This yields in particular that 
   $ |z'|^2  \left| m_{KE}  ( - \tilde \delta, \tau' , {\xi \over |\xi|})   \right|$
    is uniformly bounded for $\tau' \in \mathbb{R}.$
   Therefore, we obtain that
   $$   |G^r_{1}(t,x)| \lesssim   \int_{| \xi | \leq \delta}  e^{ - \tilde \delta | \xi | t} |\xi|^3   \int_{|\tau' | \leq \eps_{1}/|\xi|} 
    d\tau' \, d\xi \lesssim { 1 \over t^{d+2}}.$$
     To estimate the $L^1$ norm, we argue as in the proof of Lemma~\ref{lemGr2}, 
      writing
      \begin{multline}
   \label{G1rLP} G_{1}^r = \sum_{q \leq 0} G^r_{1,q},  \qquad% \quad
    %G^r_{1, q}(t,x)= \int_{\mathbb{R}^d}  e^{ i x \cdot \xi}J_{\xi, 1}(t) \chi \left( { \xi \over \delta} \right)   \phi \left( { \xi \over 2^q} \right)\, d \xi,  
  %  J_{\xi, 1}(t)= 
   % \int_{ | \tau | \leq \eps_{1}} e^{-  \tilde \delta |\xi| t}  e^{i \tau t}  {z^2 \over z^2+1}   { m_{KE} \over (z^2+ 1 - m_{KE})}(- \tilde \delta | \xi|, \tau, \xi)\, d\tau.
 %   \end{multline}
   %  Changing $\xi$  for  $\xi/ 2^q$, we get that
%\begin{multline}
%\label{G2-G2}
 G^r_{1, q}(t,x) = 2^{qd} \mathcal{G}^r_{1,q} (T, X), \\
  \mathcal{G}^r_{1,q} (T, X)
  = \int_{\mathbb{R}^d}  e^{i X \cdot \xi} J_{\xi, 1,q} (T)   \chi \left( { 2^q  \xi \over \delta} \right)   \phi \left(  \xi \right)\, d \xi, 
  \quad  T= 2^qt, \, X= 2^q x,\\
   J_{\xi, 1, q}(T)=
    \int_{ | \tau |  \leq \eps_{1} }e^{-  \tilde \delta |\xi| T}  e^{i \tau T \over 2^q}  {z^2 \over z^2+1}   { m_{KE} \over (z^2+ 1 - m_{KE})}(
    - \tilde \delta 2^q | \xi |,  \tau,  
    2^q  \xi )\, d\tau,
    \end{multline} 
    where  the integral in $\xi$ is supported on the annulus $ 1/4 \leq | \xi | \leq 4.$
    We use again  that since $2^q \leq 1$, we have the estimate 
    $$ \| \mathcal{G}_{1,q} (T, \cdot)\|_{L^1}
     \lesssim  \sum_{ | \alpha | \leq d+1}   \| \partial_{\xi}^\alpha J_{\xi, 1, q}\|_{L^\infty}.$$
     From the same estimates as above, we obtain 
       that for $1/4 \leq | \xi | \leq 4$, 
       $$    | \partial_{\xi}^\alpha J_{\xi, 1, q}(T) | \lesssim  \int_{ | \tau |  \leq \eps_{1} } |z|^2 e^{-  \tilde \delta |\xi| T} 
         \sum_{k \leq | \alpha |} 2^{qk} | h_{k} ( - \tilde \delta 2^q | \xi |,  \tau,  
    2^q  \xi )  |\, d\tau$$
    where $h_{k}$ is positively homogeneous of degree $-k$. 
%    Moreover, by using the form of $m_{KE}$ provided by 
%    Lemma \ref{remarkmKE}, we also have that 
%     for bounded $ \xi$, 
%   $$| h_{k}(z, \xi) | \lesssim { 1 \over 1 + \tau^2} .$$
%         Therefore, by setting $\tau= 2^q | \xi|  \tau'$, we obtain
%     that
Arguing exactly as in the proof of Lemma~\ref{lemGr2}, we get
   \begin{equation*}  | \partial_{\xi}^\alpha J_{\xi, 2, q}(T) | \lesssim  2^{3 q} \int_{  \tau' | \leq \eps_{1}/( 2^q | \xi|) } e^{-  \tilde \delta |\xi| T} 
         \sum_{k \leq | \alpha |} ( 1 + (\tau')^2) | h_{k} ( - \tilde \delta,  \tau',  
     {\xi \over |\xi| } )  |\, d\tau'  \lesssim   2^{2 q}  e^{- \tilde \delta T/2}.  
     % \lesssim  2^{2q}  e^{-  \tilde \delta T/2}.
      \end{equation*}
      We can then conclude as in the proof of the previous lemma by summing over the dyadic blocks.
The proof is complete.  
    \end{proof} 
    
   It  remains to estimate $G_{3}^r.$
   \begin{lem}
   \label{lemGr3}
     Under the assumptions of Proposition \ref{propGr},  $\delta$ can be chosen small enough so that  uniformly for $t \geq 1$
    $$ \| G^r_{3} \|_{L^\infty} \lesssim {1 \over t^{d+1}}, \quad   \| G^r_{3} \|_{L^1} \lesssim {1 \over t }.$$
   \end{lem}
   \begin{proof}
We are now integrating on $| |\tau |- 1 | \leq 1/2$. We shall decompose $J_{\xi,3}$ from \eqref{splitIxi1} as
   $$ J_{\xi, 3} = J_{\xi, 3, 0} + \sum_{ 1\leq k \leq N} J_{\xi, 3, k}$$
   where 
 $$ J_{\xi, 3, 0} =  \int_{ || \tau|- 1 | \leq \eps_{3} | \xi|}  e^{- \tilde \delta | \xi | t }  e^{i \tau t }
 {z^2 \over z^2+1} { m_{KE} \over (z^2+ 1 - m_{KE})}(z, \xi)
  \, d\tau$$
  and for $1 \leq k\leq N$, 
 $$ J_{\xi, 3, k} = \int_{  2^{k} \eps_{3} |\xi| \leq | | \tau|- 1 | \leq 2^{k+1} \eps_{3} | \xi|}  e^{- \tilde \delta | \xi | t }  e^{i \tau t }
 {z^2 \over z^2+1} { m_{KE} \over (z^2+ 1 - m_{KE})}(z, \xi)
  \, d\tau,$$
  where $\eps_{3}$ is given by v) of Proposition \ref{lem-zeroes}, and $N$ is such that
   $ 2^{N+1} \eps_{3}| \xi| \leq  {1 \over 2}.$

   Let us first estimate $ J_{\xi, 3, 0} $. 
    We shall focus on the estimates close to $\tau= 1$ (and call the corresponding term $J_{\xi, 3, 0}^+ $), the ones close to $\tau=-1$ can be obtained from the same arguments.
     We set
     $$ J_{\xi, 3, 0}^+  =   \int_{ | \tau- 1 | \leq \eps_{3} | \xi|}  e^{- \tilde \delta | \xi | t }  e^{i \tau t }
      {z^2 \over z^2+1} { m_{KE} \over (z^2+ 1 - m_{KE})}(z, \xi) \, d\tau
     , \quad z= - \tilde \delta |\xi|  + i \tau.$$
     We first set $\tau= 1 + | \xi | \tau'$ so that 
     $$ J_{\xi, 3, 0}^+ =  | \xi|  \int_{|\tau' | \leq \eps_{3}} e^{ | \xi | \mathfrak{z} t }  e^{i t} { ( i + | \xi | \mathfrak{z})^2 \over   2 i | \xi | \mathfrak{z} + |\xi|^2 \mathfrak{z}^2}
      { m_{KE} ( i + |\xi| \mathfrak{z}, \xi) \over | \xi| f_+ (\mathfrak{z}, |\xi|, \omega)}\, d\tau',
      $$
      where we have set $\mathfrak{z}=  (z- i)/|\xi| = -\tilde \delta + i \tau'$ and $\omega = \xi / |\xi|$ so that $\tau'= {\operatorname{Im}}\, \mathfrak{z}$ 
      and  $f_{+}$ is defined in \eqref{f+def}. %Note that from Lemma \ref{remarkmKE}, we can also write 
As in~\eqref{defm2}, we can write
       $$ m_{KE} ( i + |\xi| \mathfrak{z}, \xi)  ={ r^2  \over (i + r \mathfrak{z})^2}(3 H_{\mu}\omega \cdot \omega  +  r^2 m_{2}(i+ r \mathfrak{z}, r, \omega)), \quad r= |\xi|, $$
       therefore, we have 
       $$ |m_{KE} ( i + |\xi| \mathfrak{z}, \xi)| \lesssim |\xi|^2.$$
      Moreover,  since 
       $$ f_+ (\mathfrak{z}, r, \omega)= 2 i \mathfrak{z} + r \mathfrak{z}^2 -
        { r\over ( i + r\mathfrak{z} )^2} ( 3 H_{\mu}\omega \cdot \omega +  r^2  m_{2}(i+ r\mathfrak{z}, r, \omega)),$$
        we observe that for  $\operatorname{Re} \, \mathfrak{z} = - \tilde \delta= - \delta^{ 3 \over 2}$ and $r \leq \delta$ (see again the proof of Lemma~\ref{lemGL}), we have for
         $\delta$ small enough that
      $$| f_+ (\mathfrak{z}, r, \omega)| \gtrsim 1.$$
   By using also that for $\delta$ sufficiently small, 
   $$ \left| 2 i | \xi | \mathfrak{z} + |\xi|^2 \mathfrak{z}^2 \right| \gtrsim  |\xi|,$$
      this yields,
   $$
     |J_{\xi, 3, 0}^+| \lesssim | \xi | e^{ -\tilde \delta | \xi| t}.
    $$
The same arguments apply for $J_{\xi, 3, 0}^-$. We thus have  
    \begin{equation}
    \label{J3+}
      |J_{\xi, 3, 0}| \lesssim | \xi |e^{ -\tilde \delta | \xi| t} .
  \end{equation}   
  
   Let us now estimate $J_{\xi, 3, k}$.  Thanks to iv) of Proposition \ref{lem-zeroes}, we have
     \begin{equation}
     \label{below4} |z^2 + 1- m_{KE}| \gtrsim  2^k \eps_{3}  | \xi|.
     \end{equation}
    As above,  the estimate 
        $$ |m_{KE}| \lesssim | \xi|^2,$$
     still holds, and since $  { 1 \over 2} \geq | |\tau| - 1 | \geq 2^k \eps_{3}|\xi|$, we also have
     $$
         \quad { 1 \over |z^2 + 1|}
     \lesssim { 1 \over  2^k \eps_{3} | \xi|}, $$
   therefore,    we obtain  that
  \begin{equation}
  \label{J3k}  |J_{\xi, 3, k}| \lesssim e^{ - \tilde \delta | \xi| t} { 2^k \over  2^{2k }} |\xi | \lesssim e^{ - \tilde \delta | \xi| t} | \xi|
   { 1 \over 2^k}.
   \end{equation}
     By combining, \eqref{J3+} and \eqref{J3k}, we thus  obtain
   $$ |J_{\xi, 3}(t,x) | \lesssim  e^{ - \tilde \delta | \xi| t} | \xi| \sum_{k \geq 0} {1 \over 2^k} \lesssim e^{ - \tilde \delta | \xi| t} | \xi|.$$
   From the definition of $G^r_{3}$ (see \eqref{splitGr}, \eqref{splitIxi1}), we finally obtain
   $$ | G^r_{3}(t,x) | \lesssim { 1 \over t^{d+1}}.$$ 
   
   To estimate the $L^1$ norm, we argue again as in the proof of Lemma~\ref{lemGr2}, 
   %we use again the homogeneous Littlewood-Paley decomposition as before, 
writing
      \begin{multline}
   \label{G1rLP-3} G_{3}^r = \sum_{q \leq 0} G^r_{3,q}, \quad
   % G^r_{3, q}(t,x)= \int_{\mathbb{R}^d}  e^{ i x \cdot \xi}J_{\xi, 3} \chi \left( { \xi \over \delta} \right)   \phi \left( { \xi \over 2^q} \right)\, d \xi,  \\
  %  J_{\xi, 3}(t)= 
   % \int_{| | \tau | - 1 |\leq 1/2} e^{- | \tilde \delta |\xi| t}  e^{i \tau t}  {z^2 \over z^2+1}   { m_{KE} \over (z^2+ 1 - m_{KE})}(- \tilde \delta | \xi|, \tau, \xi)\, d\tau.
  %  \end{multline}
   %  Changing $\xi$  for  $\xi/ 2^q$, we get that
%\begin{multline}
%\label{G2-G2}
 G^r_{3, q}(t,x) = 2^{qd} \mathcal{G}^r_{3,q} (T, X), \\
  \mathcal{G}^r_{3,q} (T, X)
  = \int_{\mathbb{R}^d}  e^{i X \cdot \xi} J_{\xi, 3} (T)   \chi \left( { 2^q  \xi \over \delta} \right)   \phi \left(  \xi \right)\, d \xi, 
  \quad  T= 2^qt, \, X= 2^q x,\\
   J_{\xi, 3, q}(T)=
    \int_{ || \tau |- 1 |  \leq 1/2 }e^{-  \tilde \delta |\xi| T}  e^{i \tau T \over 2^q}  {z^2 \over z^2+1}   { m_{KE} \over (z^2+ 1 - m_{KE})}(
    - \tilde \delta 2^q | \xi |,  \tau,  
    2^q  \xi )\, d\tau,
    \end{multline} 
    where  the integral in $\xi$ is in  the annulus $ 1/4 \leq | \xi | \leq 4.$  
     To estimate  $J_{\xi, 3, q}(T)$, we focus again on the vicinity of $1$ and call the corresponding contribution
$J_{\xi, 3, q}^+(T)$. We now use the same decomposition as before for the estimate of the $L^\infty$ norm, which yields 
\begin{multline*}
 J_{\xi, 3, q}^+(T) = \sum_{0 \leq k \leq N} J_{\xi, 3, q, k}^+(T)  \\=  \sum_{0 \leq k \leq N} \int_{ 2^{k+q} |\xi| \eps_{3} \leq  |\tau - 1 | \leq \eps_{3} |\xi| 2^{k + 1+q}} 
      e^{-  \tilde \delta |\xi| T}  e^{i \tau T \over 2^q}  {z^2 \over z^2+1}   { m_{KE} \over (z^2+ 1 - m_{KE})}\left(
    - \tilde \delta 2^q | \xi |,  \tau,  
    2^q  \xi \right)\, d\tau.
    \end{multline*}
    Let us estimate
  $  \|\partial^\alpha_{\xi} J_{\xi, 3, q}^+(T)\|_{L^\infty}$ for $| \alpha | \leq d+1.$
  
   By~\eqref{defm2}, we  have that  for $1/4 \leq | \xi | \leq 4$, 
   $2^{k+q}  |\xi| \eps_{3} \leq  |\tau - 1 | \leq \eps_{3}  |\xi| 2^{k + 1+q}$ or  $ |\tau - 1 | \leq \eps_{3} 2^{2+q}$ for $k=0$,
   $$  m_{KE} (
    - \tilde \delta 2^q | \xi |,  \tau,  
    2^q  \xi )= { 2^{2q} \over   (- \tilde \delta 2^q| \xi| + i \tau)^2} \left( 3 H_\mu \xi \cdot \xi  + | \xi| ^4 m_{2} (
    - \tilde \delta 2^q | \xi |,  \tau,  
    2^q  \xi ) \right)$$
    and hence for $|\tau - 1 |\leq 1/2$ and $|\xi|\leq \delta$, we get 
        $$ \left|\partial_{\xi}^\beta \left(   {m_{KE} \over g} (
    - \tilde \delta 2^q | \xi |,  \tau,  
    2^q  \xi ) \right)\right| \lesssim  \sum_{\ell \leq | \beta | + 1 }   { 2^{2q +q(\ell-1)} \over \left|g   \left(
    - \tilde \delta 2^q | \xi |,  \tau,  
    2^q  \xi  \right)\right|^\ell},
   $$
    where we have set
    $$ g  (
    - \tilde \delta 2^q | \xi |,  \tau,  
    2^q  \xi )= (z^2+ 1 -m_{KE}) (
    - \tilde \delta 2^q | \xi |,  \tau,  
    2^q  \xi ).$$
  Consequently, by using  \eqref{below4} which gives for $|\xi |\geq 1/4$ %after the change of variable in $\xi$
 $$  | g(- \tilde \delta 2^q | \xi |,  \tau,  
    2^q  \xi ) | \gtrsim  2^{k+ q}, $$
   we obtain that
  $$  \left|\partial_{\xi}^\beta \left(   {m_{KE} \over g} (
    - \tilde \delta 2^q | \xi |,  \tau,  
    2^q  \xi ) \right)\right| \lesssim { 2^q \over 2^k}.$$
     In a similar way, we have uniformly in $q$,
  $$\left | \partial_{\xi}^\beta \left(   {(-  \tilde \delta  2^q| \xi| + i\tau)^2 \over  (-  \tilde \delta  2^q| \xi| + i\tau)^2 + 1 } \right)
   \right| \lesssim { 1 \over 2^{k+q}}.$$
   Therefore, we obtain that for $| \alpha | \leq d+1.$
 $$    \|\partial^\alpha_{\xi} J_{\xi, 3, q}^+(T)\|_{L^\infty} \lesssim  {2^{k+q} \over  2^{2k}}  e^{-  \tilde \delta T/2}
  \leq { 2^q \over 2^k}  e^{-  \tilde \delta T/2}.
  $$
  The same estimate holds for $J_{\xi, 3, q}^-(T)$. By summing over $k \geq 0$, we get that
  $$ \| \mathcal{G}_{3, q}^r  (T ) \|_{L^1} \lesssim  2^q e^{ -\tilde \delta  T/2}$$
  and we finally obtain the claimed estimate of $ \| G^r_{3} (t)  \|_{L^1}$ by summing over $q \leq 0$ as in the proof of Lemma~\ref{lemGr2}.
  
   \end{proof}
 \subsection*{End of the proof of Proposition \ref{propGr}}
 It suffices to recall the expression \eqref{splitGr} and to gather the estimates of Lemma \ref{lemGr1}, Lemma \ref{lemGr2}
  and Lemma \ref{lemGr3} (taking $\delta$ small enough).
 \end{proof}
 
%    \begin{rem}The terms $G_2^r$ and $G_3^r$ have exactly the same decay as in the screened case studied in \cite{HNR}.
% We observe however  an improved decay, namely a gain of a factor $1/t$, for $G_1^r$ (which we recall corresponds to the low frequency regime in $(\tau,\xi)$).  
%  \end{rem}

  \subsubsection{Low frequency estimates: singular part}

 We shall now study $G^S_{\pm}$ defined in \eqref{GLdec}, which corresponds to the dispersive part.
 
 \begin{prop}
 \label{propGS}
 Assuming \eqref{eq:norma}, \eqref{muanal}, \eqref{muodd} and  \eqref{muradial},  $\delta$ can be chosen small enough so that
 \begin{equation}
 \label{estGS} \|G^S_{\pm}(t) \|_{L^2}\leq C, \quad \| G^S_{\pm} (t) \|_{L^\infty} \leq {C \over t^{d \over 2}}, \quad \forall t \geq 1
 \end{equation}
 where $C$ depends on at most $ d+  1  $ derivatives of the amplitude $a_\pm$ and $ d+ 2  $ derivatives of the phase $Z_\pm$.
   We also have the more precise structure 
   \begin{equation}
   \label{plusprecis} G^S_{\pm}(t)*_{x} \cdot = e^{\pm i t} H^S_{\pm}(t, D)
   \end{equation}
   where  the operator $H^S_{\pm}(t,D)$ is such that for every $k \geq 0$, 
   \begin{equation}
   \label{higherpm}  \partial_{t}^k H^S_{\pm} (t, D)= H^S_{\pm, k}(t, D) \chi(D) \Delta^k, 
   \end{equation}
   and  $H^S_{\pm, k}(t, D)$ also satisfy the estimates
  \begin{equation}
  \label{estHS}
  \|H^S_{\pm, k}(t, D) \|_{L^2\rightarrow L^2}\leq C, \quad \| H^S_{\pm, k} (t, D) \|_{L^1 \rightarrow L^\infty} \leq {C \over t^{d \over 2} }, \quad \forall t \geq 1.
  \end{equation}
  
 \end{prop}
 
 \begin{proof}
 We focus on the study of $G^S_{+}$, the analysis of $G^S_{-}$ being similar. The estimate for the $L^2$ norm is just a consequence of the fact that the inverse
 Fourier transform is an isometry.
   We recall 
   $$ G^S_{+}(t,x) = \int_{\mathbb{R}^d} e^{Z_{+}(r, \omega) t + i x \cdot \xi} a_{+}(r,\omega) \chi\left({\xi \over \delta} \right)
    d \xi. $$
    Since we assume that~\eqref{muradial} holds, we can use Lemma \ref{lemsmooth}, from which we deduce that $Z_{+}$ is 
     a smooth function of the $\xi$ variable in $B(0, 10 \delta)$  so that we actually have 
     $$  G^S_{+}(t,x) = \int_{\mathbb{R}^d} e^{Z_{+}(\xi) t + i x \cdot \xi} a_{+}(\xi) \chi\left({\xi \over \delta} \right)
    d \xi$$
    where the  the amplitude
    $$ a_{+}(\xi)=  { Z_{+}(\xi)^2 \over  2 Z_{+}(\xi) - \partial_{z}m_{KE}(Z_{+}(\xi), \xi)}$$
    is also a smooth function of $\xi$.

   To get the decay estimate in $L^\infty$, we shall use that the imaginary part  of $Z_{+}$ described in Lemma  \ref{lemsmooth}
     provides dispersive properties. 
    Since we have almost no information on the real part of $Z_{+}$ (besides the fact that it is non-negative), we shall use a robust version of the stationary
    phase.   By using Lemma   \ref{lemsmooth}, we can write
    $$   G^S_{+}(t,x)=: e^{it} H^S_{+}(t,x)= e^{it} I(t, X), \quad X= x/t$$
    where
    \begin{equation}
    \label{defHS+} H^S_{+}(t,x)=  I(t,X) = \int_{\mathbb{R}^d}  e^{ i t  \Psi_{X}(\xi) }  a_{+}(\xi) \chi\left({\xi \over \delta} \right)
    d \xi
    \end{equation}
    with the phase given by 
    $$ \Psi_{X}(\xi) = |\xi|^2 \Phi_{+}(\xi) +   X \cdot \xi= \Psi_{X}^r(\xi) + i \Psi_{X}^i(\xi).$$
    Note  that  $ \Psi_{X}^i  \ge 0$ and 
    $$ D_{\xi}^2\Psi^r_{X} (0)=   2  {C_{\mu}} \mbox{I}_{d}.$$
  We can take  $\delta $ small enough so that
  \begin{equation}
  \label{phasebelow0}   D_{\xi}^2\Psi^r_{X} (\xi) \geq    {1 \over 2}  {C_{\mu}} \mbox{I}_{d} \geq c_{0}>0
  \end{equation}
  for $| \xi| \leq  10 \delta$
  and hence  that for every $\xi_{1}, \, \xi_{2}  \in B(0, 10\delta)$, 
    \begin{equation}
    \label{phasebelow} | \nabla  \Psi^r_{X}(\xi_{1}) - \nabla  \Psi^r_{X}(\xi_{2}) | \geq  \tilde{ c_{0}}| \xi_{1} - \xi_{2}|,  
    \end{equation}
    where the lower bound is independent of $X$.
    We will rely on the approach of Lemma 3.1 of \cite{Farah-Rousset-Tzvetkov} by checking that the imaginary part
     is harmless.
 We use the operator
$$
\mathrm{L}(u)=\frac{1}{i(1+t|\nabla \Psi_{X}|^2)}\sum_{j=1}^{d} \partial_{j}\overline{\Psi_{X}}\partial_{j}u
+
\frac{1}{(1+t |\nabla \Psi_{X}|^2)}u 
$$
(where $| \cdot |$ denotes in this context the hermitian norm of $\mathbb{C}^d$),
which satisfies by construction 
\begin{equation}\label{id-L}
\mathrm{L}(e^{it\Psi_{X}})=e^{it\Psi_{X}}
\end{equation}
and has a formal adjoint $\widetilde{\mathrm{L}}$ (i.e. $\int_{\mathbb{R}^d}\mathrm{L}u v= \int_{\mathbb{R}^d} u \widetilde{\mathrm{L}} v $, $\forall u,v \in \mathscr{C}^\infty_{c}$)  given by
\begin{multline*}
\widetilde{\mathrm{L}}(u)=-\sum_{j=1}^{d}\frac{\partial_{j}\overline{\Psi_X}}{i(1+t |\nabla\Psi_{X}|^2)}\partial_{j}u+
\Big(
-\sum_{j=1}^{d}\frac{\partial_{j}^{2}\overline{\Psi_X}}{i(1+ t |\nabla\Psi_{X}|^2)}+ 
\sum_{j=1}^{d}\frac{2t\partial_{j}\overline{\Psi_X}\,{\operatorname{Re} }(\nabla\Psi_{X}\cdot\nabla\partial_{j} \overline{\Psi_{X}})}{i(1+t|\nabla\Psi_{X}|^2)^2}\Big)u
\\+\frac{1}{(1+t |\nabla\Psi_{X}|^2)}u.
\end{multline*}
%with the notation $\nabla \phi \cdot= \sum_{j} \partial_{j} \phi\,  \partial_{j}$ and 

Using \eqref{id-L} repeatedly, we thus get  that
$$ |I(t, X) | \lesssim \int_{\mathbb{R}^d} \left|(\widetilde{\mathrm{L}})^N \left( a_{+} (\cdot)\chi\left(\frac{\cdot}{\delta}\right)\right) \right|\, d\xi,$$
for any integer $N\ge 1$. 
 We can then  check that we  get as in the proof of  Lemma 3.1 in  \cite{Farah-Rousset-Tzvetkov}
  that  
  $$(\widetilde{\mathrm{L}})^N = \sum_{ | \alpha | \leq N} a_{\alpha}^{(N)} \partial^\alpha$$
   where the coefficients $a_{\alpha}^{(N)}$ satisfy on the support of the amplitude  the estimate
   $$ |  a_{\alpha}^{(N)}| \leq C( \Lambda_{N+1}) { 1  \over   \langle t^{ 1\over 2}  \nabla \Psi_{X} \rangle^{N}}$$ with
   $$ \Lambda_{k}= \sup_{\xi \in B(0, 5 \delta)} \sup_{ 2 \leq | \alpha | \leq k} | \partial^\alpha \Psi_{X}|.$$
   Note that since $\Lambda_{k}$ involves only derivatives of order larger than $2$ of $\Psi_{X}$, this quantity
    is independent of $X$. 
    Then, by choosing $N = d+ 1 $, we get
    $$  |I(t, X) | \lesssim  C( \Lambda_{N+1}, A_{N}) \int_{B(0, \delta)}   { 1  \over   \langle t^{ 1\over 2}  \nabla \Psi_{X} \rangle^N} \, d\xi$$
    with $A_{N} = \sup_{|\alpha | \leq N} \|\partial^\alpha a_{+}\|_{L^\infty(B(0, \delta)}.$
     To conclude, we just use that
     $$    \int_{B(0, \delta)}   { 1  \over   \langle t^{ 1\over 2}  \nabla \Psi_{X} \rangle^{N}} \, d\xi
      \leq \int_{B(0, \delta)}   { 1  \over   \langle t^{ 1\over 2}  \nabla \Psi_{X}^r \rangle^{N}} \, d\xi.
$$
We finally observe that by \eqref{phasebelow0}, \eqref{phasebelow}, the map $\xi \mapsto \nabla \Psi^r_{X}$ is a diffeomorphism on $B(0, \delta)$
 and we can thus use the change of variables $\eta = \nabla\Psi^r_{X}$ and apply the bound from below of the Jacobian provided
 by \eqref{phasebelow0} to get
 $$    |I(t, X) | \lesssim  C( \Lambda_{N+1}, A_{N}) \int_{\mathbb{R}^d} { 1 \over (1 + t |\eta|^2)^{N/2} }
 \, d \eta  \lesssim   C( \Lambda_{N+1}, A_{N}) { 1 \over t^{d\over2}}.$$
This yields \eqref{estGS}. 

    To get \eqref{higherpm}, it suffices to notice that each time we take a time derivative of $H^S_{\pm}$ (see
     \eqref{defHS+}), we multiply the amplitude by $i |\xi|^2 \Psi_{+}(\xi)$. Since $\Psi_{+}$ is smooth the new amplitude
    has the same properties as before. The expression \eqref{plusprecis} follows by switching from kernels to operators. The proof is finally complete.
     \end{proof}

 \section{Proof of Theorem \ref{theokernel}}
We use \eqref{Gsplit1},  take $\delta$ small enough  
so that to apply Lemma \ref{lemGL} and the estimates of Proposition \ref{propGr} and Proposition \ref{propGS}, and finally  apply Proposition \ref{lemGH}.
Theorem \ref{theokernel} follows, with the ``regular'' part of the kernel given by
  $$ G^{R}= G^H + G^r.$$
  
\section{Proof of Theorem~\ref{thm2}}
 From the  method of characteristics,  if $f(t,x,v)$ solves \eqref{linVP}, then 
 $\rho(t,x)=\int_{\R^d} f(t,x,v)\, dv$ solves  \eqref{def-m} and hence \eqref{eqrho}
 with
 $$ S(t,x) = \int_{\R^d} f_{0}(x-vt, v) \, dv.$$
 We have the well-known dispersive estimates (see e.g. \cite{BD})
 \begin{equation}
 \label{sourceestimates}
\begin{aligned}
\| S(t)  \|_{L^1} \lesssim \| f_0 \|_{ L^1_{x,v}}, \quad  \| S(t)  \|_{L^\infty} \lesssim \frac{1}{t^d} \| f_0 \|_{ L^1_{x} L^\infty_{v}}, \quad t \geq 1, \\
 \| \na S(t)  \|_{L^1} \lesssim  \frac{1}{t} \| \nabla _x  f_0 \|_{ L^1_{x,v}}, \quad  \| \na S(t)  \|_{L^\infty} \lesssim \frac{1}{t^{d+1}} \| \nabla_x f_0 \|_{ L^1_{x} L^\infty_{v}}, \quad t \geq 1.
\end{aligned}
\end{equation}

   We then decompose
   $$ \rho(t,x)= \rho^{R}(t,x) + \rho^S_{+}(t,x)+ \rho^S_{-}(t,x)$$
   where $\rho^{R}(t,x)$ (resp. $\rho^S_{\pm}$ )
    solves
   $$ \rho^{R}= S + G^{R} *_{t,x}  S, \qquad \rho^{S}_{\pm}=  G^{S}_{\pm} *_{t,x} S. $$
    Note from Theorem \ref{theokernel} that $G^R$ verifies the same estimates as the kernel  of the linearized screened Vlasov-Poisson system
    (see Theorem 2.1 in \cite{HNR}). Therefore, we obtain the same result as in Corollary 2.1 in \cite{HNR}:
    $$ \|\rho^{R}(t) \|_{L^1} + t^d   \|\rho^{R}(t) \|_{L^\infty} \lesssim \log(1+t)\left(\| f_0 \|_{ L^1_{x,v}} + \| f_0 \|_{ L^1_x L^\infty_{v}}\right), \quad \forall t \geq 1, $$
    and we shall thus focus on the  singular part $\rho^S_{\pm}.$ We analyse the $+$ case, the other one being similar.

        The basic estimate consists in writing, thanks to \eqref{estGS} in Proposition~\ref{propGS},
        $$ \| \rho_+^S (t) \|_{L^\infty} \lesssim  \int_{0}^{t \over 2} { 1 \over (t-s)^{d \over 2}}
         \|S(s) \|_{L^1} \, ds + \int_{t\over 2}^t \|S(s) \|_{L^2} \, ds
          \lesssim { 1 \over t^{{d \over 2}- 1}}, \, \forall t \geq 1 ,$$
which does not decay for $d=1,2$. 
%          if $d \geq 3$. If $d=2$, the right hand-side diverges with a logarithmic factor; if $d=1$, this grows like $t^{1 \over 2}$.
     Assuming additionally that  $ \langle v\rangle  \nabla_{x} f_{0} \in L^1_{x,v},  \,  \langle v\rangle\nabla_{x}  f_{0} \in L^1_{x} L^\infty_{v}$, we can improve this estimate by using the refined formula \eqref{plusprecis}.
      We have 
      $$ G^S_{+}*_{t,x} S= \int_{0}^t e^{i (t-s)} H^S_{+}(t-s, D) S(s)\, ds.$$
      By using that $ { i}\partial_{s}e^{i(t-s)} = e^{i (t-s)}$, we can integrate by parts in time to get
      \begin{multline*}
      G^S_{+}*_{t,x}  S= i  (H^S_{+}(0, D) S(t) - H^S_{+}(t, D)S(0)) \\
     +    i   \int_{0}^t e^{i (t-s)} \partial_{t}H^S_{+}(t-s, D) S(s)\, ds - i \int_{0}^t e^{i(t-s)}H^S_{+}(t-s, D) \partial_{s}S(s)\, ds.
       \end{multline*}    
%       Iterating the process, we find that
%     $$  G^s_{+}*_{t,x}= \Sigma_{1} + \Sigma_{2}$$
%     where 
%     \begin{align*}
%     &  \Sigma_{1}= \sum_{j=0}^{k-1}\left(  c_{j,k} \partial_{t}^j H^S_{+}(0, D)  \partial_{t}^{k-1-j} S(t)
%      + d_{j,k}  \partial_{t}^j H^S_{+}(t,D)  \partial_{t}^{k-1-j} S(0)\right), \\
%      & \Sigma_{2} =  \sum_{j=0}^{k}  e_{j,k}   \int_{0}^t e^{i (t-s)} \partial_{t}^jH^S_{+}(t-s, D) \partial_{t}^{k-j}S(s)\, ds 
%     \end{align*}
%     where $c_{j,k}$, $d_{j,k},$ $e_{j,k}$ are harmless complex numbers.
     To estimate 
     $$\Sigma_{1}(t) =  i  (H^S_{+}(0, D) S(t) - H^S_{+}(t, D)S(0)), $$ we can use \eqref{estHS}
      and \eqref{sourceestimates}.
       This yields
       $$ \| \Sigma_{1}(t) \|_{L^\infty} \lesssim { 1 \over  t^{{d \over 2}}}  \|  f_{0}\|_{L^1_{x,v}}.$$
        For 
        $$\Sigma_{2}(t)=  i   \int_{0}^t e^{i (t-s)} \partial_{t}H^S_{+}(t-s, D) S(s)\, ds - i \int_{0}^t e^{i(t-s)}H^S_{+}(t-s, D) \partial_{s}S(s)\, ds,$$ 
        we can rely on \eqref{higherpm} (since $\chi(D)$ is a  Fourier multiplier with compactly supported
      symbol, we shall  use that   $\chi(D) \Delta$ can be bounded by $ \chi(D) |\nabla |$). This entails
        \begin{equation*} \| \Sigma_{2}(t) \|_{L^\infty} \lesssim 
        \int_{0}^{t \over 2} { 1 \over (t-s)^{d \over 2}} ( \| \nabla S(s)\|_{L^1} + \|\partial_{t} S \|_{L^1})\, ds 
       + \int_{t \over 2}^t ( \| \nabla S(s)\|_{L^2} + \| \partial_{t} S(s)\|_{L^2} )
       \, ds .% \lesssim { 1 \over t^{d \over 2}}  \int_{0}^t { 1 \over \langle s \rangle}\, ds
%           +  \int_{t\over 2}^t { 1 \over s^{ {d \over 2} + 1}}\, ds \lesssim { \log t \over t^{d \over 2} }.
           \end{equation*}            
           We  therefore need to study decay estimates for $\pa_t S$ in $L^1$ and $L^2$.  To this end, observe that
    $$ \partial_{t} S(t,x) =  - \nabla \cdot J(t,x)$$ with 
    $$J(t,x) = \int_{\R^d}v  f^l(t,x,v) \, dv$$
    where $f^l(t,x,v)$ solves the free transport equation
      $$\partial_{t} f^l + v \cdot \nabla_{x} f^l = 0$$
      with initial data $f_{0}$.
      %Therefore if 
    %  $ \langle v\rangle  f_{0} \in L^1_{x,v},  \,  \langle v\rangle f_{0} \in L^1_{x} L^\infty_{v}$ we get
  %    that
    %  $$ \|J(t)  \|_{L^\infty} \lesssim  { 1 \over t^d}$$
    %  and more generally that
 As for~\eqref{sourceestimates}, we get that
      $$    \| \na J(t)  \|_{L^\infty} \lesssim  { 1 \over t}  \|  \langle v\rangle \nabla_{v}  f_{0}\|_{L^1_{x,v}}, \quad \| \na J(t)  \|_{L^\infty} \lesssim  { 1 \over t^{d + 1 }}  \|  \langle v\rangle \nabla_{v}  f_{0}\|_{L^1_{x}L^\infty_v}. $$
     % if  $ \langle v\rangle \partial_{x}^\alpha f_{0} \in L^1_{x,v},  \,  \langle v\rangle \partial_{x}^\alpha  f_{0} \in L^1_{x} L^\infty_{v}$.
      % Since 
      % $$   \| \partial^\alpha  J(t)  \|_{L^1} \lesssim  { 1 \over t^{| \alpha |}}, $$
     %  we also have by interpolation that
    %    $$   \| \partial^\alpha J(t)  \|_{L^2} \lesssim  {1 \over t^{ {d \over 2}+ |\alpha|}}.$$
               We thus get by interpolation that
      \begin{align*}
      \label{sourceestimates}  \| \partial_{t} S(t)\|_{L^\infty} \lesssim { 1 \over t^{d + 1 } }\|  \langle v\rangle \nabla_{v}  f_{0}\|_{L^1_{x,v}} , \quad    \| \partial_{t} S(t)\|_{L^1}
       \lesssim { 1\over t } \|  \langle v\rangle \nabla_{v}  f_{0}\|_{L^1_{x}} , \\
 \| \partial_{t}  S(t)\|_{L^2} \lesssim { 1 \over t^{ {d\over 2} + 1} } \left(  \|\langle v\rangle \nabla_{v} f_{0} \|_{L^1_{x,v}} + \|  \langle v\rangle \nabla_{v}  f_{0}\|_{L^1_{x}L^\infty_v} \right).
       \end{align*}
    %   if  $  \langle v\rangle \nabla_{x} f_{0} \in L^1_{x,v},  \,  \langle v\rangle \nabla_{x}  f_{0} \in L^1_{x} L^\infty_{v}$.
     Using also~\eqref{sourceestimates}, we deduce
           \begin{align*} \| \Sigma_{2}(t) \|_{L^\infty}  &\lesssim \left({ 1 \over t^{d \over 2}}  \int_{0}^{t/2} { 1 \over \langle s \rangle}\, ds
           +  \int_{t\over 2}^t { 1 \over s^{ {d \over 2} + 1}}\, ds\right)  \left(   \|\langle v\rangle \nabla_{v} f_{0} \|_{L^1_{x,v}} + \|  \langle v\rangle \nabla_{v}  f_{0}\|_{L^1_{x} L^\infty_{v}} \right) 
            \\
           &\lesssim { \log (1 + t) \over t^{d \over 2} } \left(   \|\langle v\rangle \nabla_{v} f_{0} \|_{L^1_{x,v}} + \|  \langle v\rangle \nabla_{v}  f_{0}\|_{L^1_{x} L^\infty_{v}} \right) 
.
           \end{align*}

           This finally yields
        $$ \|\rho^S_{+}(t) \|_{L^\infty} \lesssim { \log (1+t) \over t^{d \over 2} } \left(   \|\langle v\rangle \nabla_{v} f_{0} \|_{L^1_{x,v}} + \|  \langle v\rangle \nabla_{v}  f_{0}\|_{L^1_{x} L^\infty_{v}} \right) 
$$
and the proof of Theorem~\ref{thm2} is complete.

\section{Appendix: Radial decreasing equilibria  satisfy the stability assumption (H2)}
\label{appendix}

 In this section we shall prove
 \begin{prop}
 Let $\mu$ satisfy \eqref{eq:norma} and \eqref{muanal}. If  $\mu (v) = F\left( {  |v|^2 \over 2} \right)$ with  $F'(s) <0$, $\forall  s \geq  0$, then (H2) is verified.
 \end{prop}
  
 \begin{proof}
 We study the function $ \tau\mapsto 1-m_{KE}(i \tau, \eta)$ for $\eta \in \mathbb{S}^{d-1}.$
  By using  \eqref{mKErefined}, we get that $m_{KE}(i\tau, \eta) \rightarrow 0$ when 
   $| \tau |$ tends to $+\infty$, so that it suffices to study $ \tau\mapsto 1-m_{KE}(i \tau, \eta)$
    for bounded $\tau$.  
  We have for $\gamma >0$, $|\eta|= 1$, 
  \begin{multline*} m_{KE}(z, \eta)=   -   \int_{0}^{+ \infty} e^{-(\gamma + i \tau)s}\, i \eta \cdot   \sum_{k,l} \eta_k \eta_l \mathcal{F}_v( v_k v_l \na_v \mu)(\eta s) \, d s,
  \\= - i \int_{0}^{+ \infty} \int_{\mathbb{R}^d}  e^{ -( \gamma + i \tau  + i \eta \cdot v) t} (\eta \cdot v)^3
   F'\left({|v|^2 \over 2}\right) \, dv dt.
  \end{multline*}
  We then write $v= u \eta + w$ with $w \in \eta^{\perp}= H_{\eta}$ so that
  $$  m_{KE}(z, \eta)= - i \int_{0}^{+ \infty} \int_{\mathbb{R}}  e^{ -( \gamma + i \tau  + i  u) t} 
   u^3 \Phi'\left( {  u^2 \over 2} \right)\, du dt$$
   where
    \begin{equation}
    \label{Phidefpen}\Phi(s)= \int_{H_{\eta}} F \left(  s  + {|w|^2 \over 2} \right) \, dw.
    \end{equation}
    This yields
  $$  m_{KE}(z, \eta) = - \int_{\mathbb{R}} { \tau + u \over \gamma^2 + (\tau + u)^2} u^3 \Phi'\left( {  u^2 \over 2} \right) \, du 
   - i \gamma \int_{\mathbb{R}}  {  u^3 \over \gamma^2 + (\tau + u)^2}  \Phi'\left( {  u^2 \over 2} \right) \, du.$$
   Taking the limit $\gamma \rightarrow 0$ (following e.g. \cite[Proof of Prop. 2.1]{MV}), we get that
   $$ m_{KE}(i\tau, \eta)
   = - \mbox{p.v.} \int_{\mathbb{R}} {u^3\Phi'\left( {  u^2 \over 2} \right)  \over \tau + u} \, du - i\pi \tau^3 \Phi'\left( {  \tau^2 \over 2} \right).$$
   We then observe that for bounded $\tau$ the imaginary part vanishes only for $\tau=0$
    and in this case the real part
    is equal to
    $$ - \int_{\mathbb{R}} u^2 \Phi' \left( {u^2 + |w|^2  \over 2}\right) \, du= \int_{\mathbb{R}} \int_{\R^{d-1}} F\left( {u^2 + |w|^2  \over 2}\right)\, dw du= \int_{\mathbb{R}^d}
     \mu  dv= 1.$$
     Therefore $1 - m_{KE}(i\tau, \eta)$ vanishes only for $\tau = 0$.
     
     \bigskip
     
      Let us now  compute $\partial_{z}m_{KE}(0, \eta)$ and $\partial_{z}^2 m_{KE}(0, \eta)$.
      Following the same lines, we  first get that
    \begin{multline*}
    \partial_{z}m_{KE}(z, \eta) =  i \int_{0}^{+ \infty} \int_{\mathbb{R}}  t  e^{ -( \gamma + i \tau  + i  u) t} 
   u^3 \Phi'\left( {  u^2 \over 2} \right)\, du dt  \\=  \int_{0}^{ + \infty} \int_{\mathbb{R}}e^{ -( \gamma + i \tau  + i  u) t}  \partial_{u}\left(  u^3 \Phi'\left( {  u^2 \over 2} \right)
    \right)\, du dt
        \end{multline*}
    and therefore
 \begin{multline*}
   \partial_{z}m_{KE}(z, \eta)= - i   \int_{\mathbb{R}} { \tau + u \over \gamma^2 + (\tau + u)^2}  \partial_{u}\left(  u^3 \Phi'\left( {  u^2 \over 2} \right)
    \right)
   du  \\+  \gamma \int_{\mathbb{R}}  { 1 \over \gamma^2 + (\tau + u)^2}  \partial_{u}\left(  u^3 \Phi'\left( {  u^2 \over 2} \right)
    \right)
   \, du.  
    \end{multline*}
    Taking the limit $\gamma\rightarrow 0$ as before, we get 
  $$
   \partial_{z}m_{KE}(0, \eta)  = 
     - i  \int_{\mathbb{R}}  { 1  \over u }  \partial_{u}\left(  u^3 \Phi'\left( {  u^2 \over 2} \right)
    \right) \, du.    
   $$
   By integrating by parts as before, this yields
 $$   \partial_{z}m_{KE}(0, \eta)  =   i \int_{\mathbb{R}} u \Phi'\left( {  u^2 \over 2} \right)\, du = 0
    .$$
 
  Finally,  for  $\partial_{z}^2m_{KE}(0, \eta)$, we have
 \begin{multline*} \partial_{z}^2m_{KE}(z, \eta) = -  \int_{0}^{ + \infty} \int_{\mathbb{R}}e^{ -( \gamma + i \tau  + i  u) t}  t  \partial_{u}\left(  u^3 \Phi' \left( {  u^2 \over 2} \right)
    \right)\, du dt \\ = i   \int_{0}^{ + \infty} \int_{\mathbb{R}}e^{ -( \gamma + i \tau  + i  u) t}  t  \partial_{u}^2\left(  u^3 \Phi'\left( {  u^2 \over 2} \right)   \right)\, du dt.
    \end{multline*} 
  This yields as before
  $$  \partial_{z}^2m_{KE}(0, \eta) = -  \int_{\mathbb{R}} { 1 \over u} \partial_{u}^2\left( u^3 \Phi' \left( {  u^2 \over 2} \right)    \right)\, du= \int_{\mathbb{R}}  { 1 \over u^2} \partial_{u}\left( u^3 \Phi' \left( {  u^2 \over 2} \right) \right) du
     =  2 \int_{\mathbb{R}} \Phi' \left( {  u^2 \over 2} \right) du.$$
     By using the definition \eqref{Phidefpen}, we thus get that
  $$  \partial_{z}^2m_{KE}(0, \eta) =  2 \int_{\mathbb{R}^d} F'\left( {  |v|^2 \over 2} \right) dv \neq 0$$
 and the proof of the proposition is complete.
 \end{proof}

\bigskip

 \noindent {\bf Acknowledgements:} 
DHK was partially supported by the grant ANR-19-CE40-0004. TN was supported in part by the NSF under grant DMS-1764119, an AMS Centennial fellowship, and a Simons fellowship. FR was partially supported by the ANR  projects  ANR-18-CE40-0027
 and ANR-18-CE40-0020-01.

\bibliography{unscreened.bib}
\bibliographystyle{abbrv}

\end{document}